\documentclass{amsart}

\usepackage{amssymb}
\usepackage[all]{xy}
\usepackage{hyperref}

\usepackage{enumitem}   

\setlist[enumerate]{itemsep=.2em,topsep=.2em,leftmargin=1.25em,itemindent=2.0em}

\newtheorem{thm}{Theorem}

\newtheorem{lem}[thm]{Lemma}
\newtheorem{cor}[thm]{Corollary}

\newtheorem{prop}[thm]{Proposition}

\theoremstyle{definition}
\newtheorem{defn}[thm]{Definition}

\newtheorem{say}[thm]{}
\newtheorem{exmp}[thm]{Example}
\newtheorem{exmps}[thm]{Examples}

\newtheorem{ques}[thm]{Question}    

\newtheorem{rem}[thm]{Remark}
\newtheorem{rems}[thm]{Remarks}  

\newtheorem*{ack}{Acknowledgments}      
\newtheorem{notation}[thm]{Notation}   
  
\newtheorem{defn-thm}[thm]{Definition--Theorem}  
\newtheorem{defn-lem}[thm]{Definition--Lemma}  

\newtheorem{task}[thm]{Task}

\theoremstyle{remark}

\renewcommand{\c}[0]{{\mathbb C}}  

\renewcommand{\o}[0]{{\mathcal O}} 
\newcommand{\z}[0]{{\mathbb Z}}
\newcommand{\n}[0]{{\mathbb N}}
\renewcommand{\r}[0]{{\mathbb R}} 

\renewcommand{\a}[0]{\hat{\mathbb A}} 

\newcommand{\dd}[0]{{\mathbb D}}

\newcommand{\p}[0]{{\mathbb P}}

\newcommand{\q}[0]{{\mathbb Q}}
\newcommand{\map}[0]{\dasharrow}
\newcommand{\qtq}[1]{\quad\mbox{#1}\quad}
\newcommand{\spec}[0]{\operatorname{Spec}}

\newcommand{\gal}[0]{\operatorname{Gal}}

\newcommand{\rank}[0]{\operatorname{rank}}

\newcommand{\supp}[0]{\operatorname{Supp}}    
\newcommand{\red}[0]{\operatorname{red}}

\newcommand{\coker}[0]{\operatorname{coker}}

\newcommand{\sing}[0]{\operatorname{Sing}}    
\newcommand{\ex}[0]{\operatorname{Ex}}

\newcommand{\univ}[0]{\operatorname{Univ}}

\newcommand{\rup}[1]{\lceil{#1}\rceil}
\newcommand{\rdown}[1]{\lfloor{#1}\rfloor}

\newcommand{\simq}[0]{\sim_{\q}}

\newcommand{\tsum}[0]{\textstyle{\sum}}

\def\into{\DOTSB\lhook\joinrel\to}

\def\loccoh#1.#2.#3.#4.{H^{#1}_{#2}(#3,#4)}

\DeclareMathAlphabet{\mathchanc}{OT1}{pzc}%
                                {m}{it}

\usepackage[all]{xy}\xyoption{dvips}

\newcommand{\tprod}[0]{\textstyle{\prod}}

\newcommand{\dx}[0]{\partial_x}
\newcommand{\dy}[0]{\partial_y}
\newcommand{\qg}[0]{\ensuremath{\mathrm{qG}}}

\newcommand{\defsp}[0]{\operatorname{Def}}

\begin{document}
\bibliographystyle{amsalpha}

\title{KSB smoothings of surface pairs}
        \author{J\'anos Koll\'ar}

        \begin{abstract}  We describe KSB smoothings of log canonical surface  pairs  $(S, D)$, where $D$ is a reduced curve. In sharp contrast with the $D=\emptyset$ case, cyclic quotient  pairs  always have
          KSB smoothings, usually forming many irreducible components.
        \end{abstract}

        \maketitle

        \tableofcontents

At the boundary of the moduli space of smooth surfaces, we find surfaces with log canonical singularities. Quotient singularities form the largest---and usually most troublesome---subclass.

 The universal deformation space $\defsp(S)$ of a quotient singularity is quite complicated. The irreducible components were enumerated in \cite{ksb, MR1129040}; their number grows exponentially with the multiplicity. 
However, for the moduli of surfaces, we are interested only in those deformations of $S$ for which $K_S$ lifts to a $\q$-Cartier divisor; see Paragraph~\ref{ksb.defs.say.0}.
These are now called {\it KSB deformations,} replacing the name  {\it  \qg-deformations} used in \cite{ksb}. 

By \cite{wah-loo, ksb},  quotient singularities either have no KSB smoothings, or 
the  KSB smoothings form a single,  easy to describe, irreducible component of the deformation space;  see Paragraph~\ref{ksb.defs.say}.

The aim of this note is to study the analogous question for  the moduli theory of  simple normal crossing   surface pairs. 
The normal singularities at the boundary are log canonical pairs $(S, D)$, where $S$ itself has quotient singularities whenever $D\neq \emptyset$. If the singularity is a cyclic quotient, then $D$ can have 1 or 2 local branches, and if the singularity is a dihedral quotient, then $D$  has 1 local branch; 
see  Notation~\ref{conormal.notation} for details.
The 3  types behave quite differently.

\section{Irreducible components of KSB deformation spaces}

We enumerate the irreducible components of $\defsp_{\rm KSB}(S,D)$---the {\it space of KSB deformations}---where we now require $K_S+D$ to lift to a $\q$-Cartier divisor; see Definitions~\ref{ksb.defs.pair.defn}--\ref{ksb.defs.pair.dd}.

We use the dual graph---as in
Notation~\ref{conormal.notation}---to specify a pair  $(S, D)$.

\begin{thm} \label{cyc.nplt.enum.i} The irreducible components of
  $$
  \defsp_{\rm KSB}\bigl(\bullet \ - \ c_1  \ - \ c_2  \ - \ \cdots  \ - \  c_s  \ - \ \bullet\bigr)
  $$
  are in one-to-one correspondence with the 
irreducible components of
  $$
  \defsp\bigl(c_1  \ - \ c_2  \ - \ \cdots  \ - \  c_s\bigr),
  $$
  and, in each of them, a general deformation is  smooth.
\end{thm}

\begin{thm} \label{dih.nplt.enum.i}  The irreducible components of
  $$
  \defsp_{\rm KSB}\left(
  \begin{array}{ccc}
    & 2 &\\
     & | &\\
    2  \ - \!\!\!\!   &c_1& \!\!\!\!   - \ c_2  \ - \ \cdots  \ - \  c_s  \ - \ \bullet
  \end{array}
  \right)
  $$
  are in one-to-one correspondence with the 
irreducible components of
  $$
\defsp\bigl(c_2  \ - \ \cdots  \ - \  c_s\bigr),
\eqno{\mbox{(note that  $c_1$ is omitted)}}
  $$
   and, in each of them, a general deformation has 2 singularities,  which are   $(2 \ - \ \bullet)$. 
\end{thm}

\begin{prop} \label{cyc.plt.rig.i}
  The singularities $( c_1  \ - \ c_2  \ - \ \cdots  \ - \  c_s\ - \ \bullet)$ are KSB rigid.
  \end{prop}

\begin{rems}{\ }
\begin{enumerate}
\item  The dimensions of the irreducible components of $\defsp_{\rm KSB}(S, D)$ are computed in (\ref{main.cyclic.thm}) for the cyclic cases, and in (\ref{main.cyclic.thm.d}) for the dihedral  cases.
\item  For the cyclic quotients  in Theorem~\ref{cyc.nplt.enum.i},
  the image of
  $$
  \defsp_{\rm KSB}\bigl(\bullet \ - \ c_1 \ - \ \cdots  \ - \  c_s  \ - \ \bullet\bigr)\ \to\ \defsp\bigl(c_1  \ - \cdots  \ - \  c_s\bigr)
  $$
  is nowhere dense, with a few exceptions; see (\ref{dim.drop.cyc.say}).
\item For the dihedral pairs $(S, D)$ in Theorem~\ref{dih.nplt.enum.i},
  only some irreducible components of $\defsp(S)$ contain an
  irreducible component of $\defsp_{\rm KSB}(S, D)$, see (\ref{dih.P.mod.all}). 
  \item For quotient singularities $S$, \cite[3.14]{ksb} gives an algorithm to
enumerate the irreducible components of $\defsp(S)$, though in practice this can be very cumbersome if there are many exceptional curves. For cyclic quotients $S$ as in Theorem~\ref{cyc.nplt.enum.i}, \cite{MR1129040} gives a  better method.  As a consequence, \cite{MR1129040} shows that if $S$ has  multiplicity $m=m(S)$, then $\defsp(S)$ has at most  $\tfrac1{m-1}\tbinom{2(m-2)}{m-2}$
  irreducible components, and equality holds for `most' cyclic quotients of a given multiplicity. Note that $m(S)=2+\sum (c_i-2)$. 
    \end{enumerate}
\end{rems}

{\it Other results.} 
If we  look at pairs $(S, dD)$ for  $d<1$, the situation changes dramatically. For the pairs $(S, dD)$ as in Theorems~\ref{cyc.nplt.enum.i}--\ref{dih.nplt.enum.i}, there are either no  KSBA deformations (as in Definition~\ref{ksb.defs.pair.defn}) for any $d$, or a single
irreducible component for every $d$.
For almost all of the pairs $(S, dD)$ as in Proposition~\ref{cyc.plt.rig.i}, there are 
 either no  KSBA deformations for any $d$, or a single
 irreducible component for a unique value of  $d$, see (\ref{4.0.smooth.exmp}). For $d\in (\frac56,1]$ we give a complete answer   in Section~\ref{0-1.sec}.
 For $d\in [\frac12,1]$, all lc pairs $(S, dD)$ are listed   in \cite[Sec.3.3]{kk-singbook}. Using it, 
  the method should give a full answer for  $d\in [\frac12,1]$. I checked  only  some examples, but  did not find any other pairs $(S, dD)$ with a KSBA smoothing and  $\frac12< d<1$.  Examples with $d=\frac12$ are given in (\ref{1/2.2.exmp}).

\medskip

{\it Sketch of the proofs.}
Let $(S, D)$ be a pair as in
Theorems~\ref{cyc.nplt.enum.i}--\ref{dih.nplt.enum.i}, fix  an irreducible component of
  $\defsp_{\rm KSB}(S, D)$,  and let
$g:({\mathbf S}, {\mathbf D})\to  (0,\dd)$
be a general 1-parameter deformation in it. That is,  the fiber over $0$ is
 $({\mathbf S}, {\mathbf D})_0\cong (S, D)$, and $\dd$ denotes either a complex disc or the germ of a smooth curve.

 We prove in (\ref{smoothing.prop}) that the general fiber of $g$ is smooth in the cyclic case, and has at worst $A_1$ singularities in the dihedral case.  Thus, as we recall in (\ref{cmod.ksb.thm}), by \cite[3.5]{ksb} there is a 
proper, birational morphism $\pi:{\mathbf S}_P\to {\mathbf S}$ such that
 \begin{itemize}
 \item $K_{{\mathbf S}_P}$ is $\pi$-ample, and
   \item the central fiber 
 $S_P$ has only  Du~Val and $\a^2/\tfrac1{rn^2}(1, arn{-}1)$  singularities.
 \end{itemize}
 (See Notation~\ref{cyclic.basic.not}.
 The latter are called T-singularities  (\ref{ksb.defs.say}.1), and
 $S_P\to S$ is a 
   P-modification, see   Definition~\ref{P+M.defn}).
 
 Let ${\mathbf D}_P$ denote the birational transform of ${\mathbf D}$ on
 ${\mathbf S}_P$.  Then
 $K_{{\mathbf S}_P}+{\mathbf D}_P\simq \pi^*\bigl(K_{{\mathbf S}}+{\mathbf D}\bigr)$, hence 
 $g\circ\pi:({\mathbf S}_P+{\mathbf D}_P)\to \dd$
 is a KSB deformation of $S_P$ and also of   $(S_p, {\mathbf D}_P|_{S_P})$;
 we call these doubly KSB deformations  in Definition~\ref{2-KSB.defn}.
 We show in (\ref{K+D.KSB.P.thm}) that ${\mathbf D}_P|_{S_P}=D_P+E_P$, where
 $E_P\subset S_P$ is the reduced exceptional divisor.

 A well known argument (\ref{main.MP.old.thm}) then  gives  a natural morphism
$\defsp_{\rm KSB}(S_P, D_P+E_P) \to \defsp_{\rm KSB}(S,D)$, which is 
      finite and  birational  onto our irreducible component.
   
In the cyclic case of Theorem~\ref{cyc.nplt.enum.i}, a quick argument shows that every P-modification leads to an  irreducible component of 
$\defsp_{\rm KSB}(S,D)$, see (\ref{main.cyclic.thm}).

In the dihedral  case of Theorem~\ref{dih.nplt.enum.i},  the irreducible components of $\defsp(S)$ are enumerate in (\ref{dih.P.mod.all}), following \cite{MR1129040, MR1213355}.
Then in (\ref{DP.char.lem}--\ref{dih.P.mod.all.rem})  we describe which P-modification lead to an  irreducible component of 
$\defsp_{\rm KSB}(S,D)$.    Putting these together gives 
Theorem~\ref{dih.nplt.enum.i}.

Proposition~\ref{cyc.plt.rig.i}  directly follows from \cite[2.23]{k-modbook}, see (\ref{ord.of.K+C.say}.1). \qed

\medskip

{\it Infinitesimal computations.}  Usually the space of KSB deformations $\defsp_{\rm KSB}(S, D)$ is not reduced. For cyclic pairs, the method of  \cite{k-alt} is used to  determine its tangent space in Section~\ref{first.ord.sec}, but  a complete description is known only in a  few cases.

  \medskip

{\it Other applications.}
Many of the results  apply to more general rational singularities  $S$.  However,  not   every irreducible component of $\defsp(S)$
 is obtained from a P-modification.

The conjectures in \cite[Sec.6]{k-etc} ask whether every irreducible component of $\defsp(S)$ is obtained in a similar way, using a  notion of P-modification that is more general than the one in Definition~\ref{P+M.defn}.  A positive answer  has been known for quotient singularities  \cite[3.9]{ksb} and  for  points of multiplicity $\leq 4$ \cite{MR1135880, MR1094706}.  The recent papers \cite{park2022deformations, jeon2023deformations} develop a method to prove the conjecture in many new cases. In \cite{park2022deformations} this is illustrated by the
$W(p,q,r)$ series of singularities, but their method  applies more broadly.  (The $W(p,q,r)$ series was  discovered by Wahl around 1980,  see \cite{MR4270358} for a recent survey.)  In \cite{jeon2023deformations} the conjecture is proved for most weighted homogeneous singularities.

Especially for the M-modification version  as in \cite{be-ch} or (\ref{main.MP.old.thm}), the
singularities with a rational homology disc smoothing---classified by \cite{MR2843099}---may form the natural class to work with.

\medskip

{\it Higher dimensions.} For an lc pair $(X, D)$ of arbitrary dimension with reduced $D$, \cite{k-db3} shows that a flat deformation is a KSB  deformation iff the KSB condition (\ref{ksb.defs.pair.defn})  is satisfied by a  general surface section.
  Thus our results give  information not only for surface pairs, but for higher dimensions as well.

\begin{ack}  I thank   V.~Alexeev, K.~Altmann, P.~Hacking, H.~Hauser,  J.~Stevens and  J.~Wahl  for many   useful comments and references.
Partial  financial support    was provided  by  the NSF under grant number
DMS-1901855.
\end{ack}

\section{Quotient singularities}

The classification of surface quotient singularities and  their dual graphs are given in \cite{briesk}.
We need detailed information about the cyclic (type $A$) and dihedral (type $D$) cases; we recall these and fix our notation.
The tetrahedral, octahedral and icosahedral quotients (type $E$) will appear only in examples.

\begin{notation}
We work over the complex numbers and let $\a^2$ denote the germ of ${\mathbb A}^2$ at the origin. For our purposes, we may work with a complex analytic germ,
the spectrum of $\c[[x,y]]$, or the  spectrum of the Henselisation of
$\c[x,y]_{(x,y)}$.
For most situations in this paper one can choose  global coordinates, so we can even work with   $0\in {\mathbb A}^2$.

$S_{n,q}$ and $S^d_{n,q}$ will denote the cyclic and dihedral
quotients as in Notation~\ref{cyclic.basic.not} and \ref{dih.basic.not}.
The curves  $B_{n,q}\subset S_{n,q}$, $D_{n,q}\subset S_{n,q}$ and $D^d_{n,q}\subset S^d_{n,q}$ are defined in (\ref{conormal.notation}).

The minimal resolution of a surface $S$ is denoted by $\mu:S^{\rm m}\to S$.
For a curve $D\subset S$, its birational transform on $S^{\rm m}$ is denoted by $D^{\rm m}$.
\end{notation}

\begin{notation}[Cyclic quotients]\label{cyclic.basic.not}
  We write $S_{n,q}:=\a^2/\frac1{n}(1,q)$
  where the group action is $(x, y)\mapsto (\epsilon x, \epsilon^q y)$
  for some primitive $n$th root of unity $\epsilon$. The action is free outside the origin iff $(n, q)=1$.   In general, if $n=n_1(n,q), q=q_1(n,q)$, then
  $$
  \a^2/\tfrac1{n}(1,q)\cong \a^2/\tfrac1{n_1}(1,q_1).
  $$
  Let $q'$ denote the multiplicative inverse of $q$  modulo $n$. Then
  $S_{n,q}\cong  S_{n,q'}$;  the isomorphism interchanges the coordinates.

  For the dual graph of the minimal resolution, we use $c_i:=-(C_i^2)$, the negative of the self-intersection, to denote the vertex corresponding to the curve $C_i\subset S_{n,q}^{\rm m}$.
   For $S_{n,q}$ the dual graph is
  $$
 {\mathcal D}(S_{n,q}) =\quad  c_1\ - \ c_2 \ - \cdots \  - c_s,
  $$
    where the $c_i$ are obtained from the continued fraction expansion
    $$
    \frac{n}{q}=c_1-\cfrac{1}{c_2-\cfrac{1}{c_3-\cdots}}
    $$
    We use $[c_1, \dots, c_s]$ to denote this continued fraction, and
    $S(c_1, \dots, c_s)\cong S_{n,q}$ the resulting singularity.
      \end{notation}

\begin{say}[Resolution of cyclic quotients]  \label{res.of.cyclic.say}
An explicit construction of the minimal resolution is in \cite{jung};
  see \cite{reid-cyclic} for a very accessible treatment.
  Here we use an inductive  procedure as in \cite[2.31]{k-res}.

  For ${\mathbb A}^2_{xy}/\frac1{n}(1,q)$, the quotient of the blow-up
  of the ideal sheaf $(x^q, y)$ gives a proper, birational morphism
  $\pi_1:S_1\to  {\mathbb A}^2_{xy}/\frac1{n}(1,q)$. It is covered by 2 charts
  \begin{enumerate}
  \item (Singular chart)  ${\mathbb A}^2_{x_1y_1}/\frac1{q}(1,-n)$, where
    $x=x_1y_1^{1/n}, y=y_1^{q/n}$.
    \item (Smooth chart)  ${\mathbb A}^2_{x_2y_2}$, where
      $x=x_2^{1/n}, y=y_2x_2^{q/n}$.
  \end{enumerate}
   Iterating this blow-up gives the minimal resolution.
   The $\pi_1$-exceptional curve $E_1$ is $(y_1=0)$ (resp.\ $(x_2=0)$).
   Thus the extended dual graph is
   $$
   (x-\mbox{axis})\ -\   c_1\ - \ c_2 \ - \cdots \  - c_s\ -\ (y-\mbox{axis}).
   \eqno{(\ref{res.of.cyclic.say}.3)}
  $$
\end{say}

\begin{notation}[Dihedral quotients]\label{dih.basic.not}
  Let $S^d_{n,q}$ denote the singularity whose
  minimal resolution dual graph  is
  $$
  \begin{array}{ccc}
    & 2 &\\
     & | &\\
    2  \ - \!\!\!\!   &c_1& \!\!\!\!   - \ c_2  \ - \ \cdots  \ - \  c_s 
  \end{array}
  $$
where $\frac{n}{q}=[c_1, \dots, c_s]$.  The curves marked $2$ are denoted by $C'_0, C''_0$. 
We also use the notation
$$
S(2^2, c_1,\dots, c_s):=S^d_{n,q}.
$$
One usually assumes $s\geq 2$, though later it will be convenient to allow $s=1$.

\medskip

{\it Claim \ref{dih.basic.not}.1.}   $S^d_{n,q}$ is also the quotient of
$$
S_{N,Q}:=\a^2/\tfrac1{2q(n-q)}\bigl(1, 2n'(n-q)+1\bigr), \qtq{where} qq'=nn'+1,
$$
by the  involution  induced by  $(x, y)\mapsto (\epsilon^q y, x)$, where $\epsilon$ is any primitive  $2q(n-q)$-th root of unity.
This gives the correspondence
$$
\begin{array}{cll}
  (n,q) & \mapsto & \bigl(N=2q(n-q),\ Q=2n'(n-q)+1\bigr), \qtq{and}\\
  (N,Q) & \mapsto & \bigl(n= q+ \tfrac12(N,Q-1),\   q=N/(N,Q-1)\bigr).
\end{array}
$$
(So  $N$ is even and $Q^2\equiv 1\mod N$.)
\medskip

To see these,  contract the curves  $C'_0, C''_0$ and $C_2, \dots, C_s$, we get
$\bar S^d_{n,q}\to S^d_{n,q}$. There is 
a single exceptional curve $\cong \p^1$ with two $A_1$  points on  it, plus the singularity coming from
 $S(c_2, \dots, c_s)$. Thus $\bar S^d_{n,q}$ has a double cover, ramified only at the  $A_1$  points. The corresponding  double cover of $S^d_{n,q}$ is a cyclic   quotient singularity, 
whose dual graph is
$$
c_s \ - \ \cdots\  - \ c_2\  - \ \tilde c_1\  - \ c_2   \ - \ \cdots\  - \ c_s
\eqno{(\ref{dih.basic.not}.2)}
$$
 where $ \tilde c_1= 2(c_1-1)$.
We claim that   this is the quotient singularity $S_{N,Q}$.
To see this, start with $S_{N,Q}$ and blow up the origin.
In the $(x, \frac{y}{x})$-chart  we get
$$
{\mathbb A}^2/\tfrac1{2q(n-q)}\bigl(1, 2n'(n-q)\bigr)\cong 
{\mathbb A}^2/\tfrac1{q}\bigl(1, n'\bigr).
$$
Next, write
$\tfrac{q}{m}:=[c_2,\dots, c_s],$ so  $n=c_1-\tfrac{m}{q}$,
and note that   
$nn'\equiv -1 \mod q$, so $mn' \equiv 1 \mod q$.
Therefore
$$
\a^2/\tfrac1{q}\bigl(1, n'\bigr)\cong
\a^2/\tfrac1{q}\bigl(1, m\bigr)
$$
which gives $S(c_2,\dots, c_s)$.

The exceptional curve of the blow-up is   pointwise fixed by a subgroup of order $2(n-q)$. Taking the corresponding quotient we have  a smooth surface and the self-intersection of the curve is $-2(n-q)$. Then we have a $\frac1{q}$-action, which changes the  self-intersection to $-2(n-q)/q=-2(n/q)+2$. Then we resolve the  $\a^2/\tfrac1{q}\bigl(1, m\bigr)$ singularities. The self-intersection becomes   $-2\rup{n/q}+2=-2c_1+2$, as claimed.

Finally  $(x, y)\mapsto (\epsilon^q y, x)$  descends to the involution on 
$S_{N,Q}$. \qed 
\end{notation}

\begin{say}[KSB deformations]\label{ksb.defs.say.0}
  Let $X_0$ be a surface with quotient (more generally log canonical) singularities only. A {\it KSB deformation} of $X_0$ 
  is a flat morphism to a local scheme $g: X\to (0,B)$ such that
  $X_0\cong g^{-1}(0)$, and
  \begin{enumerate}
  \item $\omega_{X/B}^{[m]}$ is flat over $B$, and  commutes with all base changes $(0,B')\to (0, B)$ for every $m$.
  \end{enumerate}
  If $B$ is reduced, then, by \cite[2.3 and 4.7]{k-modbook} this is equivalent to 
\begin{enumerate}\setcounter{enumi}{1}
  \item $\omega_{X/B}^{[m]}$ is locally free for some  $m>0$.
\end{enumerate}
The notion was introduced in \cite{ksb}, and called {\it  $\qg$-deformation} there. Note also that the 2 versions are not equivalent already over
$B=\spec \c[\epsilon]$; see \cite{k-alt}. 
\end{say}

\begin{say}[KSB smoothings]\label{ksb.defs.say}
  By \cite[5.9]{wah-loo} and \cite[3.10]{ksb}, quotient singularities that admit a KSB smoothing---frequently called {\it T-singularities---}are either Du~Val or of the form
  $$
  \a^2_{uv}/\tfrac1{rn^2}(1, arn-1),\qtq{where $(a, n)=1$ and $n>1$.}
  \eqno{(\ref{ksb.defs.say}.1)}
  $$
  More generally, these are the singularities that have a KSB deformation with Du~Val generic fiber.
  (There are a few more cases that have  nontrivial KSB deformations, but no KSB smoothings; see \cite[2.29]{k-modbook}.)
  By \cite{MR83h:14029}, the dual resolution graphs are obtained as follows.
  We start with
  $$
   4 
   \qtq{or}
         3 \ - \ 2 \ - \ \cdots \ - \ 2\ - \ 3,
    $$
    and successively apply the operation
    $$
         c_1 \ - \ \cdots \ - \   c_{s-1}  \ - \   c_s 
  \quad  \mapsto \quad
      2  \ - \   c_1  \ - \ \cdots \ - \   c_{s-1}  \ - \   (c_s+1), 
    $$
      or its symmetric version.

      Starting with $(4)$ gives the $\a^2_{uv}/\tfrac1{n^2}(1, an-1)$
      singularities, and starting with $(3 \ - \ 2 \ - \ \cdots \ - \ 2\ - \ 3)$
      (with $r{-}2$ curves marked $2$) gives the
      $\a^2_{uv}/\tfrac1{rn^2}(1, arn-1)$ singularities.

    Taking the quotient of (\ref{ksb.defs.say}.1)
    by the subgroup of order $rn$ shows that
   $$
    \a^2_{uv}/\tfrac1{rn^2}(1, arn-1)\cong (xy=z^{rn})/\tfrac1{n}(1, -1, a),
    \eqno{(\ref{ksb.defs.say}.2)}
  $$
  where $x=u^{rn}, y=v^{rn}, z=uv$. 
  Using the second representation, the universal KSB deformation is given by
  $$
  \bigl(xy=z^{rn}+\tsum_{i=0}^{r-1}t_iz^{in}\bigr)/\tfrac1{n}(1, -1, a, {\mathbf 0}).
   \eqno{(\ref{ksb.defs.say}.3)}
   $$
   Thus the KSB deformation space is smooth and has dimension $r$.

   For $r=1$ we get a 1-dimensional KSB deformation space. 
   As observed by \cite{be-ch}, the $r>1$ cases can be reduced to the $r=1$ cases as follows.
Change the deformation (\ref{ksb.defs.say}.3) to the factored form
  $$
  \bigl(xy=\tprod_{i=1}^{r}(z^{n}-s_i)\bigr)/\tfrac1{n}(1, -1, a, {\mathbf 0}),
  $$
and  repeatedly blow up the divisors $(x=z^{n}-s_i=0)$
  to obtain a small, cerpant modification, which is a flat KSB deformation of its central fiber. The central fiber is a modification of 
  $\a^2/\tfrac1{rn^2}(1, arn{-}1)$  that has $r$ singularities of the form
$\a^2/\tfrac1{n^2}(1, an{-}1)$, connected by $r-1$ curves whose birational transforms are $(-1)$-curves on the minimal resolution of the central fiber. Using Notation~\ref{boxed.nota}, the dual graph  is the following.
$$
  {\boxed{\tfrac1{n^2}(1, an{-}1)}-}\ 1\ {-\boxed{\tfrac1{n^2}(1, an{-}1)}-}\  \cdots \ {-\boxed{\tfrac1{n^2}(1, an{-}1)}-} \ 1\ {-\boxed{\tfrac1{n^2}(1, an{-}1)}}
   \eqno{(\ref{ksb.defs.say}.4)}
   $$
   This construction leads to the notion of M-modifications in (\ref{P+M.defn}).

   Note that, as we repeatedly contract $(-1)$-curves in  (\ref{ksb.defs.say}.4), we never contract the curves on the two ends. These show the following.
   \medskip

   {\it Claim \ref{ksb.defs.say}.5.}  The minimal resolution of
   $\a^2/\tfrac1{rn^2}(1, arn{-}1)$ and the minimal resolution of (\ref{ksb.defs.say}.4) are isomorphic along the birational transforms of the $x$ and $y$ axes. \qed
   \medskip

   {\it Claim \ref{ksb.defs.say}.6.}
 The map  $\a^n_{\mathbf s}\to \a^n_{\mathbf t}$, given by the elementary symmetric functions, is a Galois cover of degree $n!$. \qed
  \end{say}

  \begin{notation} \label{boxed.nota} For $\frac{n}{q}=[c_1,\dots, c_s]$ we use
    $$
    a_1\ {-\boxed{\tfrac{n}{q}}-}\ a_2   \qtq{resp.}
    a_1\ {-\boxed{\tfrac{n}{q}}}
    $$
    to denote 2 (resp.\ 1) curves through a singular point of type $\a^2/\frac1{n}(1,q)$,  whose resolution dual graph is
    $$
    a_1\ -\ c_1 \ -\  \ \cdots \ - \ c_s \ -\ a_2
    \qtq{resp.}
     a_1\ -\ c_1 \ -\  \ \cdots \ - \ c_s.
    $$
    \end{notation}

\section{Surface pairs with reduced boundary}

\begin{notation} \label{conormal.notation}
  Let  $(S, D)$ be an lc surface pair, where $D$ is a reduced curve.
 There are 3 possible local normal forms at points  $p\in D\subset S$; see for example \cite[Sec.3.3]{kk-singbook}. 
\begin{enumerate}
\item   {\it (Cyclic, plt)}
  $(S_{n,q}, B_{n,q}):=\bigl(\a^2_{xy}, (x=0)\bigr)/\tfrac1{n}(1,q)$, where $(n,q)=1$.  The minimal resolution dual graph is the following, where $\frac{n}{q}=[c_1,\dots, c_s]$ and $\bullet$ denotes the birational transform of $B_{n,q}$. 
  $$
  c_1  \ - \ c_2  \ - \ \cdots  \ - \  c_s \ - \ \bullet
     $$
\item
   {\it (Cyclic, non-plt)}
   $(S_{n,q}, D_{n,q}):=\bigl(\a^2_{xy}, (xy=0)\bigr)/\tfrac1{n}(1,q)$, where $(n,q)=1$.
     The minimal resolution dual graph is
   $$
  \bullet \ - \  c_1  \ - \ c_2  \ - \ \cdots  \ - \  c_s \ - \ \bullet,
  $$
  where $\bullet$ denotes a branch of the birational transform of $D_{n,q}$.
  Note that $D_{n,n{-}1}$ is the Du~Val singularity  $(uv=w^n)$. These are the only ones for which $D_{n,q}$ is Cartier and they behave exceptionally in many respects.
\item
  {\it (Dihedral)}  The pair  $(S^d_{n,q}, D^d_{n,q})$ has dual graph.
  $$
  \begin{array}{ccc}
    & 2 &\\
     & | &\\
    2  \ - \!\!\!\!   &c_1& \!\!\!\!   - \ c_2  \ - \ \cdots  \ - \  c_s  \ - \ \bullet,
  \end{array}
  $$
 where  $\bullet$ denotes the birational transform of $D^d_{n,q}$.
  The curve $D^d_{n,q}$ is irreducible and smooth.
  Note that $s=1$ iff $q=1$. The singularity $S^d_{n,1}\cong S_{4n-4, 2n-1}$ is a cyclic quotient, but the curve $D^d_{n,1}$ is different from the curves
  $B_{4n-4, 2n-1}$ and $D_{4n-4, 2n-1}$.

As in (\ref{dih.basic.not}.1), $(S^d_{n,q},D^d_{n,q})$  is also obtained as the quotient of
$$
(S_{N,Q}, D_{N,Q}):=\bigl(\a^2, (xy=0)\bigr)/\tfrac1{2q(n-q)}\bigl(1, 2n'(n-q)+1\bigr)
$$
by the  involution  induced by  $(x, y)\mapsto (\epsilon^q y, x)$.
 \end{enumerate}
\end{notation}

\begin{defn}[KSBA deformations of pairs]\label{ksb.defs.pair.defn}
  Let $(X_0, D_0)$ be a normal, log canonical surface pair. A {\it KSB deformation} of $(X_0, D_0)$ 
  is a flat morphism to a local scheme $g: X\to (0,B)$, plus a divisor $D\subset X$  that is Cartier at the smooth points of $X_0$,   such that
  \begin{enumerate}
  \item $X_0\cong g^{-1}(0)$,
    \item $D$ is flat over $B$  and $D_0:=D|_{X_0}$,
  \item $\omega_{X/B}^{[m]}(mD)$ is flat over $B$, and  commutes with all base changes $(0,B')\to (0, B)$ for every $m\in\z$.
  \end{enumerate}
  If $B$ is reduced, then, by \cite[2.3 and 4.7]{k-modbook} the latter  is equivalent to 
\begin{enumerate}\setcounter{enumi}{3}
  \item $\omega_{X/B}^{[m]}(mD)$ is locally free for some  $m>0$.
\end{enumerate}
If  $D$ is replaced by a $\q$- or  $\r$-divisor $\Delta$ and $B$ is reduced, then (\ref{ksb.defs.pair.defn}.2--4) should be changed to
\begin{enumerate}\setcounter{enumi}{4}
  \item $K_{X/B}+\Delta$ is a $\q$- or  $\r$-Cartier divisor.
\end{enumerate}
These are called {\it KSBA deformations.}
See \cite[Sec.8.2]{k-modbook} for the rather more complicated definition over arbitrary bases.
A major difficulty is that, even for pairs $(S, dD)$, the divisor $D$ need not be flat over the base.

\medskip
{\it Comment on terminology \ref{ksb.defs.pair.defn}.6.}
 I follow  \cite{k-modbook} for KSB and KSBA, but the distinction is important mainly over nonreduced bases, and can be safely ignored for the current purposes.
The usage in the literature  is not consistent; this applies to  my papers as well.
\end{defn}

\begin{say}[Deformation spaces]\label{ksb.defs.pair.dd}
  Let $Y$ be an affine variety (or Stein space) with isolated singularities and
  $p:X\to Y$  a proper, birational morphism such that $\ex(p)$ is proper.
  Then a universal deformation space $\defsp(X)$ exists (and is finite dimensional); see \cite{MR54:10261, bing}.
  We get $\defsp(Y)$ using the identity map $Y=Y$.
  
  If $Z\subset X$ is a closed subvariety with isolated singularities, then $\defsp(X, Z)$ exists.

  If $x\in X$ is an isolated singular  point of $X$ and $Z$, then $\defsp(x, X, Z)$ denotes the deformation space  of a small enough  neighborhood of $x\in U\subset X$ (using the identity map $U= U$).

  If $D\subset X$ is a divisor  with isolated singularities, 
all KSB deformations give a  functor. It has a universal deformation space, denoted by
$$
\defsp_{\rm KSB}(X,D).
\eqno{(\ref{ksb.defs.pair.dd}.1)}
$$
As the simplest example, adding a divisor $D$ to (\ref{ksb.defs.say}.2)
gives a  smooth irreducible component of
$\defsp_{\rm KSB}\Bigl(\bigl(\a^2_{uv}, (uv=0)\bigr)/\frac1{rn^2}(1, ran-1)\Bigr)$
$$
  \bigl((xy=z^{rn}+\tsum_{i=0}^{r-1}t_iz^{in}), (z=0) \bigr)/\tfrac1{n}(1, -1, a, {\mathbf 0}).
   \eqno{(\ref{ksb.defs.pair.dd}.2)}
   $$
   Our main aim is to describe the reduced structure of $\defsp_{\rm KSB}(X,D)$ for log canonical surface pairs.

   If  $D$ is replaced by a $\q$- or  $\r$-divisor $\Delta$, we get
   $$
\defsp_{\rm KSBA}(X,\Delta).
\eqno{(\ref{ksb.defs.pair.dd}.3)}
$$
\end{say}

As in  \cite{ksb} or \cite[2.23]{k-modbook}, there is a simplifying step if
$K_{X}+D$ is not Cartier. We study this possibility next.

\begin{say}[The order of $K_S+D$ in the class group]\label{ord.of.K+C.say}
  We consider separately the  3
  cases in  (\ref{conormal.notation}.1--3).
\begin{enumerate}
\item  For (\ref{conormal.notation}.1), the group acts faithully on 
$\frac{dx\wedge dy}{x}$. Thus  $(\a^2_{xy}, (x=0)\bigr)$ is the index 1 cover of $(S_{n,q}, B_{n,q})$; see \cite[2.49]{kk-singbook}. Therefore, as in \cite[2.23]{k-modbook},
every KSB deformation of $(S_{n,q}, B_{n,q})$ is induced by a KSB deformation of  $(\a^2_{xy}, (x=0)\bigr)$. The latter is rigid, so every
KSB deformation of $(S_{n,q}, B_{n,q})$ is trivial.
\item  For (\ref{conormal.notation}.2),
  the rational 2 form
  $\frac{dx\wedge dy}{xy}$ is invariant under the group action, hence descends to a generating section of  $\omega_{S_{n,q}}(D_{n,q})$.  Thus
  $\omega_{S_{n,q}}(D_{n,q})$ is locally free and $K_S+D$ is Cartier.
  Thus the method of \cite[2.23]{k-modbook} does not give any information about the  KSB deformations.

\item  For (\ref{conormal.notation}.3),  the rational 2 form
  $\frac{dx\wedge dy}{xy}$ is invariant under the cyclic subgroup of order $2q(n-q)$, but anti-invariant under  $(x, y)\mapsto (\epsilon^q y, x)$. Thus only its tensor square descends, so
  $K_S+D$ is 2-torsion in the local Picard group.
  The corresponding duble cover is $(S_{N,Q}, D_{N,Q})$ as in (\ref{conormal.notation}.3).

  As in \cite[2.23]{k-modbook},
  every KSB deformation of $(S^d_{n,q}, D^d_{n,q})$ is induced by a KSB deformation of  $(S_{N,Q}, D_{N,Q})$. Thus, in principle, this reduces us to the previous case. We give a more direct description  in (\ref{dih.P.mod.all}).
 
\end{enumerate}

Pulling $K_S+D$ back to the minimal resolution shows that
\begin{enumerate}\setcounter{enumi}{3}
 \item  $K_{S^{\rm m}}+D^{\rm m}+\sum_{i=1}^s C_i\sim 0$ in case (\ref{conormal.notation}.2), and
\item   $2\bigl(K_{S^{\rm m}}+D^{\rm m}+\sum_{i=1}^s C_i\bigr)+(C'_0+C''_0) \sim 0$  in case (\ref{conormal.notation}.3).  
\end{enumerate}
\end{say}

\section{Q-modifications}

P-modifications were introduced in \cite{ksb} (under the name P-resolution) to describe   smoothings of quotient singularities. A variant, using M-modifications is developed in  \cite{be-ch}.  To handle arbitrary deformations of pairs we introduce 
Q-modifications.  These are actually very natural objects from the point of view of birational geometry, and many basic results on P- and M-modifications
hold for Q-modifications. 

\begin{defn}\label{Q-mod.defn}
  Let $(s,S)$ be a rational surface singularity. A proper, birational morphism $\pi:S_Q\to S$ is a {\it Q-modification}
if
\begin{enumerate}
\item $K_{S_Q}$ is $\pi$-nef, and
\item $S_Q$ has only quotient singularities. (So `Q' for quotient.)
\end{enumerate}
\end{defn}

\begin{notation}\label{P.res.ample.say}
Let $(0,S)$ be a rational singularity with minimal resolution
$\mu:S^{\rm m}\to S$. Let $\pi:\bar S\to S$ be a proper, birational morphism.
Assume that  $\bar S$ is normal with  minimal resolution
 $\bar\mu:\bar S^{\rm m} \to \bar S$.
Then  $\bar S^{\rm m}$ dominates $S^{\rm m}$, giving a commutative diagram
  $$
    \xymatrix{%
      \bar S^{\rm m}  \ar[r]^{\bar\mu} \ar[d]_{\pi^m}   & \bar S  \ar[d]^{\pi} \\
S^{\rm m}\ar@{.>}[ur]^{\tau}  \ar[r]^{\mu}   &  S
    }
    \eqno{(\ref{P.res.ample.say}.1)}
 $$   
    We can write
    $$
    K_{\bar S^{\rm m}}\simq  \bar\mu^* K_{\bar S}-\bar\Delta,\qtq{and}
    K_{\bar S^{\rm m}}\sim (\pi^m)^*K_{S^{\rm m}}+F,
    \eqno{(\ref{P.res.ample.say}.2)}
    $$
    where $\bar\Delta$ is effective with $\supp \bar\Delta\subset \ex(\bar\mu)$,
    and $F$ is effective  with $\supp F = \ex(\pi^m)$.

        Let $\bar A^{\rm m}\subset \bar S^{\rm m}$ be an irreducible  curve that is exceptional over $S$, with  images  $A^{\rm m}\subset  S^{\rm m}$ and $\bar A\subset \bar S$. (We use $0$ if the image is  a point.)
    Then (\ref{P.res.ample.say}.2) gives that
$$
      (\bar A\cdot K_{\bar S})=
      (A^{\rm m}\cdot K_{S^{\rm m}})+\bigl(\bar A^{\rm m}\cdot F\bigr)+
      \bigl(\bar A^{\rm m}\cdot \bar \Delta\bigr).
      \eqno{(\ref{P.res.ample.say}.3)}
      $$
      From this we conclude the following.
      \begin{enumerate}\setcounter{enumi}{3}
      \item  If $\bar A\neq\emptyset$ and $ A^{\rm m}\neq\emptyset$, then
        $(\bar A\cdot K_{\bar S})\geq 0$, with equalty holding iff
        $(A^{\rm m}\cdot K_{S^{\rm m}})=0$ and $\tau$ is an isomorphism along
        $A^{\rm m}$.
      \item If $ A^{\rm m}=\emptyset$, then
        $ (\bar A\cdot K_{\bar S})$ depends only on the intersection numbers of the curves in $\ex(\bar\mu)\cup \ex(\pi^m)$.
      \end{enumerate}
\end{notation}

\begin{lem}\label{antieff.lem} Using the notation of  (\ref{P.res.ample.say}), assume that $K_{\bar S}$  is $\pi$-nef. 
  Then 
    $-K_{\bar S^{\rm m}}$ is $\q$-linearly equivalent to an effective linear combination of curves, that are exceptional for $\bar S^{\rm m}\to S$.
\end{lem}

Proof. By (\ref{P.res.ample.say}.2) $ K_{\bar S^{\rm m}}\simq  \bar\mu^* K_{\bar S}-\bar\Delta$ where $\Delta $ is effective. Since $K_{\bar S}$  is $\pi$-nef,
$-K_{\bar S}$  is $\q$-linearly equivalent to an effective linear combination of $\pi$-exceptional curves. \qed
\medskip

\begin{rem}  \label{antieff.lem.rem} This is especially strong if $S$ has quotient singularities. Then, by \cite[3.13]{ksb}, there is a unique, maximal
  $\pi':S'\to S$ such that $-K_{S'}$  is $\q$-linearly equivalent to an effective linear combination of $\pi'$-exceptional curves.
  This gives an effective procedure to construct all Q-modifications.
  \end{rem}

\begin{lem}\label{rigid.quot.mod.lem.1}
  Let $\pi:\bar S\to S$ be a Q-modification and  $\pi^m$ as in (\ref{P.res.ample.say}.1).
   Then   $\pi^m$ is a composite of a  
  sequence of blow ups
  $$
  \bar S^{\rm m}=:Y_r\to Y_{r-1}\to \cdots \to Y_1\to Y_0:= S^{\rm m},
  $$
  where each $Y_{i+1}\to Y_i$ is the blow-up of a singular point of the
  image if $\bar\Delta$ in $Y_i$.
\end{lem}

Proof. By (\ref{P.res.ample.say}.2)   $K_{\bar S^{\rm m}}\sim  \bar\mu^* K_{\bar S}-\bar\Delta$, 
 where
$\rdown{\bar\Delta}=0$ since $\bar S$ has only quotient singularities.

 For any $j$, let  $\bar\Delta_j$ denote the push-forward of $\bar\Delta$ to $Y_j$. Then   $K_{Y_r}+\bar\Delta_r\simq \bar\mu^* K_{\bar S}$ is nef over $S$
 by assumption, hence so is each
 $K_{Y_j}+\bar \Delta_j$.
 Let now $F_{i+1}$ denote the 
 exceptional curve of $Y_{i+1}\to Y_i$. Then $(K_{Y_{i+1}}\cdot F_{i+1})=-1$,
  so $(\bar \Delta_{i+1}\cdot F_{i+1})\geq 1$.  Since $\rdown{\bar\Delta_j}=0$, $F_{i+1}$ must intersect at least 2
  irreducible components of $\bar\Delta_{i+1}$ (different from $F_{i+1}$).
  The images of these in $Y_i$  show that the blow-up center is
  a singular point of $\bar\Delta_{i}$.
\qed

\medskip

These show that Q-modifications are combinatorial objects that are
determined by the dual graph ${\mathcal D}(S)$.
In fact, the following much stronger result holds, which generalizes  \cite[7.2]{MR1129040}.
To state it,
consider a graph with vertices $V$ and edges $E$. For a subset  $V'\subset V$ the  {\it induced 
  subgraph}  contains all the edges in $E$ that go  between vertices in $V'$.

\begin{prop}
  \label{Pmod.extends.prop}
  Let $S_1, S_2$ be surfaces with rational singularities
  and $j: {\mathcal D}(S_1)\into {\mathcal D}(S_2)$ an
isomorphism onto an induced 
subgraph that preserves the labeling of vertices.
Then there is a one-to-one correspondence between
    \begin{enumerate}
      \item Q-modifications $\tau_1: S_1^{\rm m}\to S_1\leftarrow  \bar S_1$ of $S_1$, and
      \item Q-modifications  $\tau_2: S_2^{\rm m}\to S_2\leftarrow \bar S_2$ of $S_2$,  such that $\ex(\tau_2)$ is  supported on
        $j( {\mathcal D}(S_1))$.
    \end{enumerate}
    Under this correspondence, $\bar S_1$ and $\bar S_2$ have the same singularities.
\end{prop}

Proof. Given a Q-modification  of $S_1$, we construct  a Q-modification  of $S_2$ as follows.

By (\ref{rigid.quot.mod.lem.1}),  $\pi_1^{\rm m}: \bar S_1^{\rm m}\to S_1^{\rm m}$ is obtained by repeatedly blowing up singular points of the exceptional divisor. We do the same blow ups to
get $\bar S_2^{\rm m}\to S_2^{\rm m}$. Then $j$ lifts to an
embedding $\bar{\jmath}$ of the dual graph of $\ex(\bar S_1^{\rm m}\to S_1)$
into the dual graph of $\ex(\bar S_2^{\rm m}\to S_2)$.

Next we get $\bar\mu_1: \bar S_1^{\rm m}\to \bar S_1$ 
by contracting certain
curves contained in
$\ex(\bar S_1^{\rm m}\to S_1)$.
Using $\bar{\jmath}$, we construct
$\bar\mu'_2: \bar S_2^{\rm m}\to \bar S'_2$
by contracting the  corresponding
curves of $\ex(\bar S_2^{\rm m}\to S_2)$.
By construction, $\bar S_2$ has the same singularities as $\bar S_1$.

We need  to show that $K_{\bar S_2}$ has nonnegative intersection number with
every exceptional curve of $\bar S_2\to S_2$.
For images of curves in $\ex(\pi_2^{\rm m})$ this follows from (\ref{P.res.ample.say}.5). For images of curves in $\ex(\mu_2)$ this follows from (\ref{P.res.ample.say}.4).

Conversely, let  $\tau_2: S_2^{\rm m}\to S_2\leftarrow \bar S_2$  be a  Q-modification of $S_2$ as in (\ref{Pmod.extends.prop}.2).
Then 
(\ref{rigid.quot.mod.lem.1}) shows that
$\bar S_2^{\rm m}\to S_2^{\rm m}$ is obtained by repeatedly blowing up
singular points of $j( {\mathcal D}(S_1))$. This shows that the above procedure can be reversed to get $\tau_1: S_1^{\rm m}\to S_1\leftarrow  \bar S_1$ of $S_1$.
\qed

\medskip

\begin{defn} \label{sim.Q.mod.defn} Let $p:X\to \dd$ be a flat morphism such that $X_0$ has rational singularities. A {\it simultaneous Q-modification} is a
  proper, birational morphism  $\pi:\bar X\to X$ such that
\begin{enumerate}
       \item $\pi_t:\bar X_t\to X_t$ is a Q-modification for every $t\in \dd$, and
    \item $K_{\bar X}$ is $\q$-Cartier. Equivalently,     $p\circ \pi:\bar X\to \dd$ is a KSB deformation.
\end{enumerate}
Being a Q-modification is an open property for KSB deformations.
Thus if  $\pi_0:\bar X_0\to X_0$ is a Q-modification, then the same holds for nearby $t$.
\end{defn}

The next result, proved in \cite[3.5]{ksb},  says that every flat deformation of a surface  with quotient singularities  has a simultaneous Q-modification.
(See  \cite[5.41]{k-modbook} for higher dimensions.)

\begin{thm}\label{cmod.ksb.thm}
Let $p:X\to \dd$ be a flat morphism such that $X_0$ has quotient singularities. Then there is a unique  simultaneous Q-modification $\pi:\bar X\to X$ such that
  \begin{enumerate}
  \item $\pi$ is an isomorphism over $X\setminus X_0$, and 
        \item $K_{\bar X}$ is  $\pi$-ample.
    \item If $X_t$ has Du~Val singularities for $t\neq 0$, then
 $\pi_0:\bar X_0\to X_0$ is a P-modification  (\ref{P+M.defn}).
      \qed
  \end{enumerate}
\end{thm}

\section{P-modifications}

By \cite[3.9]{ksb}, for every quotient singularity $S$, the irreducible components of $\defsp(S)$ are in one-to-one correspondence with the 
P-modifications of $S$. We define and study these next.

\begin{defn}\label{P+M.defn}
  Let $(s,S)$ be a rational surface singularity. A proper, birational morphism $\pi_P:S_P\to S$ is a {\it P-modification}
if $S_P$ is normal,
\begin{enumerate}
\item $K_{S_P}$ is $\pi_P$-ample, and
\item $S_P$ has only Du~Val and $\a^2/\tfrac1{rn^2}(1, arn{-}1)$ singularities, where  $(a,n)=1$.
\end{enumerate}
Following  \cite{be-ch}, $\pi_M:S_M\to S$ is an
 {\it M-modification} if  $S_M$ is normal,
\begin{enumerate}
\item $K_{S_M}$ is $\pi_M$-nef, and
\item $S_M$ has only  $\a^2/\tfrac1{n^2}(1, na{-}1)$ singularities, where  $n>1$ and $(a,n)=1$.
\end{enumerate}
Given a P-modification  $S_P\to S$, resolving its Du~Val  singularities and applying the construction (\ref{ksb.defs.say}) 
gives an M-modification $\pi_{MP}:S_M\to S_P$ such that
$K_{S_M}\simq \pi_{MP}^*K_{S_P}$. By \cite{be-ch} this establishes a
one to one correspondence between P- and M-modifications of $S$.

{\it Note on terminology.} For P-modifications we follow \cite[3.8]{ksb}.
A more general version is defined in \cite[6.2.13]{k-etc}, which allows more  singularities on $S_P$.
\end{defn}

For P-modifications  we get the following variant of (\ref{Pmod.extends.prop}).

\begin{cor}\label{Pmod.extends.prop.cor} Using the notation and assumptions of (\ref{Pmod.extends.prop}), there is a one-to-one correspondence between
    \begin{enumerate}
      \item P-modifications $\tau_1: S_1^{\rm m}\map  \bar S_1$ of $S_1$, and
      \item P-modifications  $\tau_2: S_2^{\rm m}\map \bar S_2$ of $S_2$,  such that $\ex(\tau_2)$ is  supported on
        $j( {\mathcal D}(S_1))\cup\{(-2)\mbox{ curves disjoint from }j( {\mathcal D}(S_1))\} $.
    \end{enumerate}
\end{cor}

Proof. The only difference is that, as we go from $\bar S_1$ to $\bar S_2$, we have to ensure that $K_{\bar S_2}$ is ample. That is, we have to contract
each connected component of  all $(-2)$-curves disjoint from $j( {\mathcal D}(S_1))$ to  a Du~Val singularity. \qed

\begin{say}[P-modifications of dihedral quotients]\label{dih.P.mod.all}
  By \cite{MR1129040, MR1213355}
  P-modifications $\tau_D: S_D^{\rm m}\to S_D\leftarrow \bar S_D$ of $S_D:=S(2^2, {\mathbf c})$ are obtained from
  P-modifications  $\tau_C: S_C^{\rm m}\to S_C\leftarrow \bar S_C$ of  $S_C:=S(2, {\mathbf c})$ as follows.
    \begin{enumerate}
  \item ($C_1$ is not contracted by $\tau_C$.)  Then  $C'_0, C''_0$ are  contracted to  $A_1$ points by $\tau_D$.
  \item ($C'_0,C_1, \dots, C_i$ are contracted by $\tau_C$ to an $A_{i+1}$ point.)  Then  $C''_0$ is also  contracted  by $\tau_D$, giving a  $D_{i+2}$ point.
    \item ($C_1, \dots, C_i$ are contracted by $\tau_C$ to a non-DV point.)  Then  $C'_0, C''_0$ are  not contracted   by $\tau_D$.
    \item ($C'_0, C_1, \dots, C_i$ are contracted by $\tau_C$ to a non-DV point.)  Then  $C''_0$ is  not contracted  by $\tau_D$. There is also a symmetric version, where
      $C''_0,C_1, \dots, C_i$ are contracted by $\tau_D$ but $C'_0$ is  not contracted.
    \end{enumerate}
    \end{say}

  The P-modifications that appear in Theorem~\ref{dih.nplt.enum.i}  have several characterizations.

  \begin{prop} \label{DP.char.lem}
    Let $\pi_P:S_P\to S$ be a P-modification of a dihedral quotient singularity
    $(S^d_{n,q}, D^d_{n,q})$ with reduced exceptional curve $E_P$. The following are equivalent.
    \begin{enumerate}
    \item $S_P$ is one of the  cases (\ref{dih.P.mod.all}.1--2).
      \item $C'_0, C''_0$ are contracted to  Du~Val point(s) on $S_P$.
     \item   $(S_P, D_P+E_P)$ is log canonical. 
  \item   $K_{S_P}+D_P+E_P$ is numerically $\pi_P$-trivial.
    \end{enumerate}
\end{prop}

  Proof.  The equivalence of (1) and (2) is clear.
  To see (3) note that if $C_1, \dots, C_i$ are contracted by $\tau_C$ to a non-DV point, then there are 3 curves in $D_P+E_P$ through this point, and if
  $C'_0, C_1, \dots, C_i$ are contracted by $\tau_1$ to a non-DV point, then
  $C''_0$ meets the resolution at $C_1$, which is not  an end curve  of
  the chain for $i\geq 2$. For $i=1$ it is an  end curve,  but then another curve of $D_P+E_P$ meets $C_1$.

  For (4), note that $K_{S^{\rm m}}+\tsum_{i>0} C_i+\frac12 (C'_0+C''_0)$ is numerically trivial on $S^{\rm m}\to S$. Thus
  $K_{S_P}+D_P+E_P-\frac12 (C'_0+C''_0)$ is numerically $\pi_P$-trivial.
  Thus (4) holds iff both $C'_0,C''_0$ are contracted. This happens only in 
 cases (\ref{dih.P.mod.all}.1--2).\qed

  \begin{rem} \label{dih.P.mod.all.rem} Note that the cases (\ref{DP.char.lem}) are in bijection with
    P-modifications of  $T:=S(c_2,\dots, c_s)$.
    Indeed,  let $\tau_T:  T^{\rm m}\to T\leftarrow  T_P$
  be a P-modification.

    If $c_1>2$ then we get $\tau_D$ 
    by not contracting $C_1$ but contracting  $C'_0, C''_0$. 

    If $c_1=2$ and $C_2$ is contracted to a non-DV point by $\tau_T$, then again we contract  $C'_0, C''_0$ but not $C_1$. If $C_2,\dots, C_i$ are contracted to an $A_{i-1}$  point  by $\tau_T$, then $C'_0, C''_0, C_1, \dots, C_i$ are  contracted to a $D_{i+2}$ point  by $\tau_D$.
    Finally if  $C_2$ is not contracted, then 
    $C'_0, C''_0, C_1$  are  contracted to an $A_{3}$ point  by $\tau_D$.
  \end{rem}

The main  reason to study 
P- and M-modifications is the following, essentially proved in \cite[3.9]{ksb} and  \cite{be-ch}.

\begin{thm}\label{main.MP.old.thm} Let $(0,S)$ be a  rational surface singularity,  $\pi_{P}:S^{}_P\to S$  a P-modification,  and
  $\pi_{MP}:S^{}_M\to S^{}_P$ the corresponding M-modification.
  Then there is a commutative diagram
  $$
\begin{array}{ccccc}
\univ_{\rm KSB} (S^{}_M) & \stackrel{c^{}_{MP}}{\longrightarrow} &
  \univ_{\rm KSB} (S^{}_P) & \stackrel{c^{}_{P}}{\longrightarrow} & \univ (S)\\[1ex]
 u^{}_M \downarrow\hphantom{u^{}_M}  && u^{}_P \downarrow\hphantom{u^{}_P}  && u_S \downarrow\hphantom{u_S} \\
  \defsp_{\rm KSB} (S^{}_M) & \stackrel{\tau^{}_{MP}}{\longrightarrow} &
  \defsp_{\rm KSB} (S^{}_P) & \stackrel{\tau^{}_{P}}{\longrightarrow} &
  \defsp (S)
\end{array}
\eqno{(\ref{main.MP.old.thm}.1)}
$$
where
\begin{enumerate}\setcounter{enumi}{1}
  \item $\defsp_{\rm KSB}(S^{}_M)$ and
    $\defsp_{\rm KSB}(S^{}_P)$ are smooth of the same dimension,
  \item   $\tau^{}_{MP}: \defsp_{\rm KSB}(S^{}_M) \to \defsp_{\rm KSB}(S^{}_P)$ is finite and Galois,
  \item $\tau^{}_{P}: \defsp_{\rm KSB}(S^{}_P) \to \defsp(S)$
    is a finite, birational morphism onto an irreducible component of
    $\defsp(S)$.
    \end{enumerate}
\end{thm}

Sketch of proof. For the existence of the diagram we need to show that the
morphisms  $\pi_{P}:S^{}_P\to S$   and
$\pi_{MP}:S^{}_M\to S^{}_P$ extend to
their flat deformations
$$
c^{}_{P}:  \univ_{\rm KSB} (S^{}_P)\to \univ (S)
\qtq{and}
c^{}_{MP}: \univ_{\rm KSB} (S^{}_M) \to \univ_{\rm KSB} (S^{}_P).
$$
The relevant property is the vanishing of
$R^1(\pi_{P})_*\o_{S^{}_P}$ and of $R^1(\pi_{MP})_*\o_{S^{}_M}$. 
The proof is summarized in \cite[11.4]{km-flips}, which in turn relies mainly
on \cite{MR54:10261}.

Since KSB deformations are locally unobstructed, the obstruction to smoothness lies in $H^2$ of the tangent sheaf. If $X\to U$ is proper of fiber dimension $\leq 1$, and $U$ is affine, then $H^i$ of any coherent  sheaf on $X$ vanishes for $i\geq 2$. This implies smoothness in (\ref{main.MP.old.thm}.2); see (\ref{loc.glob.lem}) for details.

The relation between
P- and M-modifications is established in \cite{be-ch}.
The Galois group in  (\ref{main.MP.old.thm}.3) is a product of reflection groups, determined by the singularities of $S^{}_P$, see \cite{be-ch} and (\ref{ksb.defs.say}.6) for a key example.

In order to show (\ref{main.MP.old.thm}.4), note first that 
every fiber of $u^{}_P$ is a P-modification of the corresponding deformation of $S$. We check in   (\ref{rigid.quot.mod.thm}) that P-modifications $\pi_{P}:S^{}_P\to S$ do not have nontrivial 1-parameter deformations fixing $S$.
This shows that  $\tau^{}_{P}$ is finite.
As in \cite{art-bri},  $\tau^{}_{P}$
    is  onto an irreducible component of $\defsp(S)$
    by openness of versality.
Thus $\tau^{}_{P}: \defsp_{\rm KSB}(S^{}_P) \to \defsp^{}(S)$
is a finite cover. In order to show that it has degree 1, let
$T\to \defsp^{}(S)$ be a morphism from the spectrum of a DVR that maps the closed (resp.\ generic) point  the closed (resp.\ generic) point.
By what we already proved, there is a finite $T'\to T$ such that the pull-back of $\univ(S)$ to $T'\to T\to \defsp^{}(S)$ is obtained from a KSB deformation of  $S^{}_P$. We may assume that $T'/T$ is Galois. As noted in \cite[p.312]{ksb}, 
the action of $\gal(T'/T)$ lifts to an action on this KSB deformation.
Taking the quotient shows that the pull-back of $\univ(S)$ to $T\to \defsp^{}(S)$ is obtained from a KSB deformation of  $S^{}_P$. Thus
$T\to \defsp^{}(S)$ lifts to $T\to \defsp_{\rm KSB}(S^{})$, so
  $\tau^{}_{P}$ is birational. \qed

\medskip
{\it Remark \ref{main.MP.old.thm}.5.} By  \cite{MR1129040},  $\tau^{}_{P}$
is  an isomorphism onto an irreducible component of $\defsp(S)$ for cyclic quotients, but
there are dihedral examples where $\tau^{}_{P}$
is  not an isomorphism.
Also, $\tau^{}_{P}$ is birational but not necessarily finite
 for the more general 
P-modifications considered in \cite[Sec.6]{k-etc}.

\section{Tangent sheaves}

In order to compute the dimensions of the deformation spaces, 
we need various results on  tangent and logarithmic tangent sheaves for
modifications of rational singularities. These are mostly taken from
\cite{MR367277, MR333238} and 
\cite{lee-park, MR2995167}.  Note that \cite{lee-park, MR2995167} focus on the global aspects, which tend to be more subtle. 
See also \cite{pa-pa-sh-2, MR4509033, MR4502069} for similar computations.

\begin{say}[Tangent and logarithmic tangent sheaves]\label{log.tg.sh.say}
  Let  $S$ be a smooth surface and $D\subset S$ a smooth curve with normal bundle $N_D$.
  There is a natural map $T_S\to N_D$. Its kernel is the
  {\it logarithmic tangent bundle} of the pair $(S, D)$, denoted by
  $T_S(-\log D)$.  Thus we have
  an exact sequence
  $$
  0\to T_S(-\log D) \to T_S\to N_{D}\to 0.
  \eqno{(\ref{log.tg.sh.say}.1)}
  $$
  The restriction of $  T_S(-\log D)$ to $D$ sits in an exact sequence
  $$
  0\to \o_D\to T_S(-\log D)|_D \to T_D\to 0.
  \eqno{(\ref{log.tg.sh.say}.2)}
  $$
  Locally, if $D=(y=0)$, then the kernel  $\o_D$ is generated by
  $y\tfrac{\partial}{\partial y}$ and the quotient $T_D$ by
  $\tfrac{\partial}{\partial x}$.

  Next assume that  $S$ is smooth but $D\subset S$ is a nodal curve
  with normalization $\tau: \bar D\to D$. The {\it immersed normal sheaf} of $D$ is
  $$
  N_{\bar D}:=  \tau_*\coker\bigl[T_{\bar D}\to \tau^*T_S\bigr].
   \eqno{(\ref{log.tg.sh.say}.3)}
 $$
 Then (\ref{log.tg.sh.say}.1) becomes
  $$
  0\to T_S(-\log D) \to T_S\to N_{\bar D}\to 0,
 $$
   where $T_S(-\log D)$ is locally free. More generally, if $D=D_1+D_2$, then we have
$$
  0\to T_S(-\log (D_1+D_2)) \to T_S(-\log D_1)\to N_{\bar D_2}\to 0.
   \eqno{(\ref{log.tg.sh.say}.4)}
 $$
   The sequence (\ref{log.tg.sh.say}.2) now gives
   $$
  0\to \o_{\bar D}\to \tau^*T_S(-\log D) \to T_{\bar D}(-N)\to 0,
  \eqno{(\ref{log.tg.sh.say}.5)}
  $$
  where $N\subset \bar D$ is the preimage of the nodes of $D$.
  
  If $S$ is a normal surface and $D\subset S$ a reduced curve, then there is a finite subset $P\subset S$ such that $(S\setminus P, D\setminus P)$ is a smooth pair. Let $j:S\setminus P\into S$ be the natural injection.
  The {\it tangent sheaf} of $S$  and the {\it logarithmic tangent sheaf} of $(S, D)$ are defined as
  $$
  T_S:=j_*T_{S\setminus P}, \qtq{and}
  T_S(-\log D):= j_*\bigl(T_{S\setminus P}(-\log (D\setminus P))\bigr).
   \eqno{(\ref{log.tg.sh.say}.6)}
   $$
   In the log canonical case, we have the following close analog of (\ref{log.tg.sh.say}.4).

\end{say}

\begin{lem} \label{log.tg.sh.thm}
  Let $(S, D_1+D_2)$ be a log canonical pair, 
   $\mu:S^{\rm m}\to S$  the minimal resolution, and $D^{\rm m}_i:=\mu^{-1}_*D_i$ with 
normalizations  $\tau_i: \bar D^{\rm m}_i\to D^{\rm m}_i$. Then 
there is an exact sequence
  $$
0\to T_S(-\log (D_1+D_2)) \to T_S(-\log (D_1))\to \mu_*N_{\bar D^{\rm m}_2}\to 0.
 \eqno{(\ref{log.tg.sh.thm}.1)}
  $$
\end{lem}

Proof.  Over the smooth locus of $S$, this is (\ref{log.tg.sh.say}.4).
Thus it remains to check what happens at the singular points that are contained in $D_1+D_2$. For this we can use local analytic coordinates.

  If a finite group $G$ acts with isolated fixed points on a surface $X$, then the local sections of $T_{X/G}$ are the $G$-invariant local sections of $T_X$.
Thus the tangent sheaf of $S_{n,q}=\a^2_{xy}/\frac1{n}(1,q)$ is generated by
  $$
  x\tfrac{\partial}{\partial x}, y^{q'}\tfrac{\partial}{\partial x}, y\tfrac{\partial}{\partial y}, x^q\tfrac{\partial}{\partial y}.
   \eqno{(\ref{log.tg.sh.thm}.2)}
   $$
   If $B_x\subset S_{n,q}$ is the image of the $x$-axis, then we see that
   $T_S(-\log B_x)$ is generated by
   $x\tfrac{\partial}{\partial x}, y^{q'}\tfrac{\partial}{\partial x},  y\tfrac{\partial}{\partial y}$, and the quotient
   $T_S/T_S(-\log B_x)$ is generated by $x^q\tfrac{\partial}{\partial y}$.

   For the minimal resolution of $S_{n,q}$ we use  the chart given in
   (\ref{res.of.cyclic.say}.2). The chain rule  gives that
   $$
   (\pi_1)_* \bigl(\tfrac{\partial}{\partial y_2}\bigr)=x^q\tfrac{\partial}{\partial y}.
   $$
   That is, $T_S/T_S(-\log B_x)$ is naturally isomorphic to the
   normal bundle of $B_x^{\rm m}\subset S^{\rm m}_{n,q}$. This shows (\ref{log.tg.sh.thm}.1) at the cyclic quotient points.

   At a dihedral point  $(S^d_{n,q}, B^d_{n,q})$, the  minimal resolution can be obtained by first taking a double cover
   $(S_{N,Q}, D_{N,Q})\to (S^d_{n,q}, B^d_{n,q})$, resolving $S_{N,Q}$ and then quotienting out by the involution. The latter is fixedpoint-free along the birational transform of $D_{N,Q}$, so the normal bundle computation for the birational transform of $B^d_{n,q}$ is the same as for the birational transform of $D_{N,Q}$.
   Finally note that $Q^2\equiv 1\mod N$, so  the 
   quotient  $T_S/T_S(-\log B_x)$ is generated by 
   $x^Q\tfrac{\partial}{\partial y}+y^Q\tfrac{\partial}{\partial x}$,
   where $x,y$ are the  orbifold  coordinates on $S_{N,Q}$.
   This shows (\ref{log.tg.sh.thm}.1) at the dihedral quotient points. \qed

   \begin{prop}\label{cyc.lc.lf.prop}
     Let $(X, D)$ be a log canonical pair with  quotient singularities.  Let $\pi:X^{\rm m}\to X$ be either the minimal resolution,
or the blow-up of a node of $D$ where $X$ is smooth.
    Let   $E^{\rm m}$ be the exceptional curve and  $D^{\rm m}$ the birational transform of $D$. Then
     \begin{enumerate}
     \item $\pi_*T_{X^{\rm m}}\bigl(-\log (E^{\rm m}+D^{\rm m})\bigr)=T_X(-\log D)$,\item $R^1\pi_*T_{X^{\rm m}}\bigl(-\log (E^{\rm m}+D^{\rm m})\bigr)=R^1\pi_*T_{X^{\rm m}}\bigl(-\log (E^{\rm m})\bigr)= 0$, and
     \item  $H^i\bigl(X, T_X(-\log D)\bigr)=
    H^i\bigl(X^{\rm m}, T_{X^{\rm m}}\bigl(-\log (E^{\rm m}+D^{\rm m})\bigr)\bigr)$ for every $i$.
     \end{enumerate}
   \end{prop}

 Proof.
 For (\ref{cyc.lc.lf.prop}.1) the 
 blow-up of a node of $D$ is a simple computation.
At a  singular point of $X$, 
 it suffices to show the $D=\emptyset$ case. That is,  all local generators of
   $T_X$  lift to sections of
 $T_{X^{\rm m}}\bigl(-\log E^{\rm m}\bigr)$.  
 This should hold for all normal singularities (\ref{TX.list.rem}), but here is an argument for quotient  singularities.

 Write $S=\a^2/G$. Let $\pi:B_0\a^2\to \a^2$ the blow-up of the origin with  exceptional curve $E$. By direct computation,
 $\pi_*T_{B_0\a^2}(-\log E)= m_0T_{\a^2}$,
 where $m_0$ is the maximal ideal at the origin.
Set $\pi':S':=B_0\a^2/G\to S$ with  exceptional curve $E'$.
 Taking $G$-invariants gives that
$\pi_*T_{S'}(-\log E')= (m_0T_{\a^2})^G$.
 Now note that the $G$-action on $T_{\a^2}/m_0T_{\a^2}$ has no fixed vectors,
 so $(m_0T_{\a^2})^G=T_{\a^2}^G=T_S$.

 The $G$-action on $B_0\a^2$ leaves $E$ invariant, so 
 $S'$ has only cyclic quotient singularities.
 Unfortunately, iterating  the above process does not give a resolution, see \cite[2.29]{k-res}. However, for cyclic quotients we can use the blow-up charts
(\ref{res.of.cyclic.say}.1--2).  By  direct computation
   $
   x\tfrac{\partial}{\partial x}=
   nx_2\tfrac{\partial}{\partial x_2}-q y_2\tfrac{\partial}{\partial y_2},
   $
   proving (\ref{cyc.lc.lf.prop}.1).

The $D=\emptyset$ case of 
(\ref{cyc.lc.lf.prop}.2)  is discussed in (\ref{taut.defn}.1).
For general $D$, we use (\ref{log.tg.sh.say}.4) to get
$$
0\to  T_{X^{\rm m}}\bigl(-\log (E^{\rm m}+D^{\rm m})\bigr)
\to  T_{X^{\rm m}}\bigl(-\log (E^{\rm m})\bigr)
\to  N_{D^{\rm m}}\to 0.
$$
Pushing this forward, we get that
$$
R^1\pi_*T_{X^{\rm m}}\bigl(-\log (E^{\rm m}+D^{\rm m})\bigr)\cong
\coker\bigl[T_S\to \pi_*N_{\bar D^{\rm m}}\bigr],
$$
  and the latter is 0 by (\ref{log.tg.sh.thm}).
   The Leray spectral sequence now gives (\ref{cyc.lc.lf.prop}.3).\qed

   \medskip
       {\it Remark \ref{cyc.lc.lf.prop}.4.}   Let $\pi:X'\to X$ be the blow-up of a smooth point with exceptional curve $E$. Then
       $h^1(X', T_{X'}(-\log E))=h^1(X, T_X)+2$.
       
   \medskip
{\it Complement \ref{cyc.lc.lf.prop}.5.}  The computations also show that if all the singularities are 
$(S_{n,q}, D_{n,q})$, then  $T_X(-\log D)$ is locally free
and $\pi^*T_X(-\log D)\cong  T_{X^{\rm m}}\bigl(-\log (E^{\rm m}+D^{\rm m})\bigr)$.

\begin{ques} \label{TX.list.rem}
  Let $X$ be a normal variety and  $\pi:X^{\rm c}\to X$ its canonical modification.  Every automorphism of $X$ lifts to an automorphism of $X^{\rm c}$, and, if $X$ is proper, then  $H^0(X, T_X)$ is the tangent space of the automorphism group.  Thus one can expect that
  $$
  \pi_*T_{X^{\rm c}}=T_X.
  \eqno{(\ref{TX.list.rem}.1)}
  $$
  \cite[4.1]{MR1320605} says that if $X$ is local and complete, then again
  $H^0(X, T_X)$ is the tangent space of the automorphism group, suggesting that 
  (\ref{TX.list.rem}.1) indeed holds.

  The caveat is that  both objects are infinite dimensional, so one would need to check that the above finite-dimensional intuition  works out in general. 
  \end{ques}

We can now compute the global contribution to the dimension of  our
deformation spaces.

\begin{prop}\label{dim.h1T.prop.ol}
  Let $S$ be a affine  surface with rational singularities and $D\subset S$ a reduced curve. Let 
  $\mu:S^{\rm m}\to S$ be the  minimal resolution with  reduced exceptional
  set $E^{\rm m}$, and 
  $\pi:{\bar S}\to S$  a Q-modification with  reduced exceptional set $\bar E$. Let  $D^{\rm m}$ resp.\ $\bar D$ denote the birational transforms of $D$ on $S^{\rm m}$  resp.\ $\bar S$. Assume that $(\bar S, \bar E+\bar D)$ is lc and
  $(S^{\rm m},E^{\rm m}+D^{\rm m})$ is normal crossing.
  
For an irreducible component $\bar E_i\subset \bar E$, let
  $e_i\in \n$ be the negative of the self-intersection of its
birational transform on the minimal resolution of ${\bar S}$.
Then
\begin{enumerate}
\item $h^1\bigl({\bar S}, T_{\bar S}(-\log \bar D)\bigr)=
  h^1\bigl({\bar S}, T_{\bar S}(-\log ( \bar E+ \bar D))\bigr)+
  \tsum_{i\in I} (e_i-1)$, and
  \item 
    $h^1\bigl({\bar S}, T_{\bar S}(-\log ( \bar E+ \bar D))\bigr)=
    h^1\bigl(S^{\rm m}, T_{S^{\rm m}}(-\log (E^{\rm m}+D^{\rm m}))\bigr).$
\end{enumerate}
\end{prop}

Proof.  
(\ref{log.tg.sh.thm}) gives an exact sequence
$$
0\to  T_{\bar S}\bigl(-\log ( \bar E+ \bar D)\bigr)\to
T_{\bar S}(-\log \bar D)\to \oplus \o_{\bar E_i}(-e_i)\to 0.
$$
Here $h^0\bigl(\bar E_i, \o_{\bar E_i}(-e_i)\bigr)=0$ since $e_i>0$, so taking cohomologies gives the first claim.
For the second, we use the notation of (\ref{P.res.ample.say}), and claim that
$$
h^1\bigl(\bar S^{\rm m}, T_{\bar S^{\rm m}}(-\log (\bar F^{\rm m}+\bar D^{\rm m}))\bigr)
$$
equals  both terms in (\ref{dim.h1T.prop.ol}.2), where $\bar F^{\rm m} $ is the reduced exceptional set of  $\bar S^{\rm m}\to   S$.
For $\bar S^{\rm m}\to  \bar S$ we use the minimal resolution case of (\ref{cyc.lc.lf.prop}).
Next, since  $\bar S^{\rm m}\to  S^{\rm m}$ is a composite of node blow-ups by (\ref{rigid.quot.mod.lem.1}), we can repeatedly use the
node blow-up case of (\ref{cyc.lc.lf.prop}).  \qed

\section{Dimension formulas}

We compute the dimension of the irreducible components of various deformation spaces using invariants of the singularity and its P-modifications.
The formulas are simplest when $S$ is determined by ${\mathcal D}(S)$.

\begin{say}
  \label{taut.defn}
Let  $(s, S)$ be a rational surface singularity with  
 minimal resolution  $\mu:S^{\rm m}\to S$ and  reduced exceptional curve
$C^{\rm m}=\cup_{i\in I} C_i$.

By \cite[Sec.3.4.4]{sernesi}, $H^1\bigl(S^{\rm m}, T_{S^{\rm m}}(-\log C^{\rm m})\bigr)$
is the tangent space to those deformations of $S^{\rm m}$ where every $C_i$ lifts. Equivalently, these are  those deformations of $S$ that preserve the dual graph.  The dimension of this tangent space is usually hard to compute, but it is known in several important cases.

\medskip
{\it Quotient singularities \ref{taut.defn}.1.} The dual graph determines
the singularity, hence $H^1\bigl(S^{\rm m}, T_{S^{\rm m}}(-\log C^{\rm m})\bigr)=0$.
This has been long known; see for example \cite{briesk},
  \cite[III.5.1]{bpv} or 
  \cite[3.32]{kk-singbook}.

  \medskip
      {\it Taut singularities \ref{taut.defn}.2.} These are the singularities that are determined by their dual graphs.
      For these $H^1\bigl(S^{\rm m}, T_{S^{\rm m}}(-\log C^{\rm m})\bigr)=0$.
      A rather laborious result of \cite{MR367277, MR333238} says that
 this vanishing implies tautness,  with a few exceptions; these are enumerated in \cite[3.2]{MR333238}.
 The   complete list of taut singularities is given in \cite{MR333238}.

 A simple case is when  ${\mathcal D}(S)$ has only 1 fork, which has  degree 3 and self-intersection $\leq -3$. This is  used for the  $W(p,q,r)$ singularities studied in \cite{park2022deformations}.

   \medskip
       {\it Weighted homogeneous singularities \ref{taut.defn}.3.}
       For these  ${\mathcal D}(S)$ has only 1 fork. Let $d_0$ denote its   degree  and  $-c_0$ its self-intersection.  By \cite[4.1.III]{MR333238}, if
       $c_0\geq 2d_0-3$ then
       $H^1\bigl(S^{\rm m}, T_{S^{\rm m}}(-\log C^{\rm m})\bigr)=d_0-3$.

        These are used
        in \cite{jeon2023deformations} under the (mostly) weaker assumption
        $c_0\geq d_0+3$.
\end{say}

For  P-modifications, we need only the simplest numerical 
invariants.

\begin{defn}\label{I.J.P.mod.defn}   Let  $\pi:S_P\to S$ be a
  P-modification of a normal singularity. Let
$$
\{p_j\in S_P: j\in J_P\} \qtq{resp.} \{E_i\subset S_P: i\in I_P\}
$$
 denote the singular points of $S_P$, respectively
 the irreducible, exceptional curves of $\pi:S_P\to S$.
 Set $r_j=r$ if $p_j$ is of type $A_r, D_r, E_r$  or 
 $\a^2/\tfrac1{rn^2}(1, arn{-}1)$ for some $n>1$.
Note that $r_j=\dim \defsp_{\rm KSB}(p_j,S_P)$ by (\ref{ksb.defs.say}.2).
 
Let $e_i\in \n$ denote the negative of the self-intersection of the
birational transform of $E_i$ on the minimal resolution of $S_P$.
\end{defn}

We can now state the first dimension formula.

\begin{prop}\label{taut.dim.thm} Let $\pi:S_P\to S$ be a  P-modification of an affine surface  with  rational  singularities, and  $I_P, J_P$  as in (\ref{I.J.P.mod.defn}).
  Let $S^{\rm m}\to S$ be the minimal resolution, with reduced exceptional curve $C^{\rm m}$. 
  Then 
  $$
  \dim \defsp_{\rm KSB}(S_P)=\tsum_{j\in J_P}r_j+\tsum_{i\in I_P} (e_i-1)+h^1\bigl(S^{\rm m}, T_{S^{\rm m}}(-\log C^{\rm m})\bigr).
  \eqno{(\ref{taut.dim.thm}.1)}
  $$
  In particular, if $S$ has  taut singularities, then
  $$
  \dim \defsp_{\rm KSB}(S_P)=\tsum_{j\in J_P}r_j+\tsum_{i\in I_P} (e_i-1).
  \eqno{(\ref{taut.dim.thm}.2)}
  $$
\end{prop}

Proof.   First note that, by (\ref{loc.glob.lem}),
$$
\dim \defsp_{\rm KSB}(S_P)= \tsum_{j\in J_P} \dim  \defsp_{\rm KSB}(p_j,S_P)
+ h^1(S_P, T_{S_P}).
$$
Here  $r_j=\dim \defsp_{\rm KSB}(p_j,S_P)$ by (\ref{ksb.defs.say}.2),
this explains the summand $\tsum_{j\in J_P}r_j$. 
 The computation of $ h^1(S_P, T_{S_P})$ is the  $D=0$ special case of  (\ref{dim.h1T.prop.ol}). \qed

\medskip
For KSB deformations of pairs we have the following.
The dihedral case is discussed in (\ref{main.cyclic.thm.d}).

   \begin{thm}\label{main.cyclic.thm} Let $\pi:(S_P, D_P+E_P)\to (S,D):=(S_{n,q}, D_{n.q})$ be a P-modification of a cyclic quotient  singularity and 
     $J_P$  as  in (\ref{I.J.P.mod.defn}).
     Then
     $$
     \dim \defsp_{\rm KSB}(S_P, D_P+E_P)= \tsum_{j\in J_P}r_j+
     |\mbox{type $A$  points}|.
     \eqno{(\ref{main.cyclic.thm}.1)}
     $$
   \end{thm}

   {\it Clarification \ref{main.cyclic.thm}.2.} Here
   $ |\mbox{type $A$  points}|$ is the number of  points 
    where the pair $(S_P, D_P+E_P)$ is locally of the form
   $\bigl((xy=z^{r+1}), (xy=0)\bigr)$ for some $r\geq 0$.
   For $r=0$ the surface is smooth but the curve is singular. 
   \medskip

   Proof. If $p_j$ is a singular point of $S_P$ or of $D_P+E_P$, then
   $\dim \defsp_{\rm KSB}(p_j, S_P, D_P+E_P)$  is  $r$ for
   a singularity is of type  $\a^2/\tfrac1{rn^2}(1, arn{-}1)$ by (\ref{ksb.defs.pair.dd}.2), but
   $r+1$ if the singularity is of type  $A_r$.  (Here we need to count
   $\bigl(\a^2, (xy=0)\bigr)$ as type $A_0$.) Thus the
   right hand side of (\ref{main.cyclic.thm}.1) is the sum of the local terms in (\ref{loc.glob.lem}).

   It remains to show that the global term
   $h^1\bigl(S_P, T_{S_P}(-\log (D_P+E_P))\bigr)$ vanishes.
   This follows from (\ref{dim.h1T.prop.ol}.2),  (\ref{cyc.lc.lf.prop}.2) and
   (\ref{taut.defn}.1.). \qed

   \begin{rem} \label{dim.drop.cyc.say} Comparing (\ref{taut.dim.thm}) and
     (\ref{main.cyclic.thm}) we see that
     $$
     \dim \defsp_{\rm KSB}(S_P, D_P+E_P)\leq \dim \defsp_{\rm KSB}(S_P),
     $$
     and the inequality is strict with a few exceptions. These are
     triple points  (\ref{KSB.codin.infini.thm}) and also P-resolutions of the form
   $$
    \bullet {-\boxed{A_{n_1}}-}\ 2\ {-\boxed{\tfrac1{rn^2}(1, arn{-}1)}-} \ 2\
     {-\boxed{A_{n_2}}-} \bullet
     $$
     \end{rem}

The dihedral version of (\ref{main.cyclic.thm}) is the following.
  
   \begin{thm}\label{main.cyclic.thm.d} Let $\pi:(S_P, D_P+E_P)\to (S,D)=(S^d_{n,q}, D^d_{n.q})$ be a P-modification of a dihedral quotient  singularity as in
     (\ref{DP.char.lem}) and 
     $J_P$  as  in (\ref{I.J.P.mod.defn}).
     Then
     $$
     \dim \defsp_{\rm KSB}(S_P, D_P+E_P)= \tsum_{j\in J_P}r_j+
     |\mbox{type $A$  points}| -2.
     \eqno{(\ref{main.cyclic.thm.d}.1)}
     $$
   \end{thm}

   Proof.  Note that here $(S_P, D_P+E_P)$ has a singular point of type
   $$
   \bigl((z^2=x(y^2-x^{r-2})), (x=z=0)\bigr)$$  
   whose KSB deformation space has dimension $r{-}2$ by  (\ref{D.term.QC.exmp.1}).  (We have 2 singular points if $r=2$.)
   This is why we need to subtract 2 on
the
   right hand side of (\ref{main.cyclic.thm.d}.1). The rest of the argument is as for (\ref{main.cyclic.thm}). \qed

\medskip
We used the  following result, whose
 proof  is summarized in \cite[11.4]{km-flips}.

\begin{lem}\label{loc.glob.lem}
  Let $S$ be an affine  surface,  
  $\tau:X\to S$ a proper, birational morphism and $D\subset X$ a reduced curve. Assume that $X$ is normal and let $\{p_j: j\in J\}$ be the   points  of $\sing X\cup \sing D$. Then the restriction maps
  $$
  \defsp (X,D)\to \times_{j\in  J} \defsp (p_j, X, D)
  \qtq{and}
  \defsp_{\rm KSB} (X,D)\to \times_{j\in  J} \defsp_{\rm KSB}  (p_j, X,D)
  $$
 are smooth of relative dimension $h^1\bigl(X, T_X(-\log D)\bigr)$.  \qed
\end{lem}

By (\ref{ord.of.K+C.say}.1), the pairs  $(S_{n,q}, B_{n,q})$ are KSB rigid. By contrast, we show that the pairs  $(S_{n,q}, D_{n,q})$ are KSB smoothable, and the pairs 
$(S^d_{n,q}, D^d_{n,q})$ have KSB deformations whose general fibers have only
singularities of type 
$(S_{2,1}, D_{2,1})\cong  \bigl((xy=z^2), (x=z=0)\bigr)$; we call these
{\it almost smoothings.} More precisely, we have the following.

\begin{prop} \label{smoothing.prop}
  Every irreducible component of
  $\defsp_{\rm KSB} (S_{n,q}, D_{n,q})$ contains smoothings, and 
every irreducible component of
   $\defsp_{\rm KSB} (S^d_{n,q}, D^d_{n,q})$ contains almost smoothings.
\end{prop}

Proof.   By (\ref{ord.of.K+C.say}.4--5),  $K_S+D$ (resp.\ $2(K_S+D)$)  is Cartier, so the same holds for its KSB deformations, hence every singularity of a KSB deformation of our cyclic (resp.\ dihedral) pair is again   cyclic (resp.\ dihedral). By openness of versality
\cite{art-bri}, it is thus enough to show that each  cyclic (resp.\ dihedral) pair is smoothable  (resp.\ almost smoothable).

We find such deformations in the Artin component.
For a cyclic pair $(S, D)$, 
 look at the minimal resolution  
  $\pi_0:(S^{\rm m}, D^{\rm m}+E^{\rm m})\to (S, D)$. By (\ref{loc.glob.lem}),
 there is a deformation  $p: ({\mathbf S}^{\rm m}, {\mathbf D}^{\rm m})\to \dd$
 whose special fiber is $(S^{\rm m}, D^{\rm m}+E^{\rm m})$ and whose generic fiber is a smooth pair with ${\mathbf D}^{\rm m}_{gen}$ irreducible.
 Next $\pi_0$ extends to a contraction
 $$
 \pi:({\mathbf S}^{\rm m}, {\mathbf D}^{\rm m})\to
 ({\mathbf S}, {\mathbf D})\to \dd,
 $$
and $({\mathbf S}, {\mathbf D})\to \dd$ is a KSB smoothing
 since $K_{S^{\rm m}}+ D^{\rm m}+E^{\rm m}\sim 0$  by (\ref{ord.of.K+C.say}.4).

 For a dihedral pair $(S, D)$,  we start with
 $\pi_0:(S', D'+E')\to (S, D)$ which is obtained from the
 minimal resolution by contracting the curves  $C'_0, C''_0$  (as in Notation~\ref{dih.basic.not}).
  The key point is that $2\bigl(K_{S'}+ D'+E'\bigr)\sim 0$
  by (\ref{ord.of.K+C.say}.5).
  The above argument now works for this case too, except that we have 2 singularities of type
 $\bigl((xy=z^2), (x=z=0)\bigr)$. These are KSB rigid, so persist in every deformation. 
\qed

\section{Doubly KSB deformations}\label{0-1.sec}

\begin{defn}\label{2-KSB.defn} A flat morphism 
  $p: {\mathbf S}\to \dd$   of $S\cong {\mathbf S}_0$  is a {\it doubly KSB deformation} of $(S, D)$ if
  it is a KSB deformation both for $S$ and for  $(S, D)$.
  Equivalently, both $K_S$ and $D$ lift to $\q$-Cartier $\z$-divisors on
  ${\mathbf S}$.

  Since we have an explicit description of all  KSB deformations,
  we can use it to describe all doubly KSB deformations in most cases. The basic examples are the following.
\begin{enumerate}
 \item Every KSB smoothing of $S_{n,q}$ is a
  doubly KSB deformation of  $(S_{n,q}, D_{n,q})$ by (\ref{ksb.defs.pair.dd}.2).
\item By Proposition~\ref{cyc.plt.rig.i}, $(S_{n,q}, B_{n,q})$ has no doubly KSB deformations.
  \item The non-Du~Val dihedral singularities  have no  KSB deformations by Paragraph~\ref{ksb.defs.say}; see (\ref{D.term.QC.exmp.1}) for the Du~Val cases.
\item 
  If  $\frac56<d\leq 1$ and $(S,dD)$ is lc, then $(S, D)$ is also lc, see \cite[3.44]{kk-singbook}.
  Thus the doubly KSB smoothings are fully described by (\ref{2-KSB.defn}.1--3).
\item Other examples with $d=\frac12$ are in (\ref{1/2.2.exmp});
  see also (\ref{D.term.QC.exmp.1}).
\end{enumerate}

\end{defn}

\begin{thm}\label{K+D.KSB.P.thm}  Let $(S, dD)$ be an lc pair such that
  either $d>\frac12 $ or $d=\frac12$ and   $(S, dD)$ is klt.
  There is a bijection between the following sets.
\begin{enumerate}
\item  \label{K+D.KSB.P.thm.1}
  Irreducible components of $\defsp_{\rm KSBA}(S, dD)$ (\ref{ksb.defs.pair.dd}.3) whose generic fiber has only Du~Val singularities.
\item \label{K+D.KSB.P.thm.2}
  P-modifications $\pi_P:S_P\to S$ such that
  \begin{enumerate}
\item  \label{K+D.KSB.P.thm.2.a} $K_{S_P}+d(D_P+E_P)\simq 0$, and
\item   \label{K+D.KSB.P.thm.2.b} the singularities of $\bigl(S_P, D_P+E_P\bigr)$ have doubly KSB deformations with Du~Val  generic fibers.
 \end{enumerate}
\end{enumerate}
\end{thm}

Proof.   Assume that  $\pi_P:S_P\to S$ satisfies (\ref{K+D.KSB.P.thm}.2).
By (\ref{loc.glob.lem})  the local doubly KSB deformations  glue together to a doubly KSB deformation of
$\bigl(S_P, D_P+E_P\bigr)$. The contraction $\pi_P$ extends to nearby deformations as in (\ref{main.MP.old.thm}), and we get a whole irreducible component by  openness of versality \cite{art-bri}.  (Note that this part holds without any restriction on $d$.)

Conversely, take a general 1-parameter KSBA deformation
$p: ({\mathbf S}, d{\mathbf D})\to \dd$ of $(S, dD)$.
Let   $\pi: \bar {\mathbf S}\to {\mathbf S}$ be the simultaneous P-modification as in (\ref{cmod.ksb.thm}.3).

Let $\bar{\mathbf D}$ be the birational transform of ${\mathbf D}$.
Since $\pi$ is small,
$$
K_{\bar {\mathbf S}}+d\bar{\mathbf D}\simq \pi^*\bigl(K_{{\mathbf S}}+d{\mathbf D}\bigr),
$$
and $(\bar {\mathbf S},d\bar{\mathbf D}+\bar {\mathbf S}_0)$ is also lc.
By  adjunction as in \cite[4.9]{kk-singbook}, we get that the pair
$\bigl(\bar {\mathbf S}_0, d\bar{\mathbf D}|_{\bar {\mathbf S}_0}\bigr)$ is also lc. Each divisor in $ d\bar{\mathbf D}|_{\bar {\mathbf S}_0}$ appears with a coefficient that is an integer multiple of $d$, and also at most $1$.
If $d>\frac12 $ then all coefficients are $d$ or $0$. If  $(S, dD)$ is klt
then all coefficients are $<1$, so again  they are  $d$ or $0$.

Since $K_{\bar {\mathbf S}_0}$ is $\pi_0$-ample,
$-\bar{\mathbf D}|_{\bar {\mathbf S}_0}$ is also $\pi_0$-ample,
hence its support must contain the whole exceptional divisor. Therefore
$\bar{\mathbf D}|_{\bar {\mathbf S}_0}$ is the sum of $\bar D$ plus the
$\pi_0$- exceptional divisor $\bar E_0$, so
$$
K_{\bar {\mathbf S}_0}+d(\bar D+\bar E_0)=\bigl(K_{\bar {\mathbf S}}+d\bar{\mathbf D}\bigr)|_{\bar {\mathbf S}_0}\simq \pi^*\bigl(K_{{\mathbf S}}+d{\mathbf D}\bigr).
 \qed
$$

 This takes us to the following, which is somewhat more general than needed for
 (\ref{K+D.KSB.P.thm}.2).

\begin{task} \label{task.task}
  Let $(S, dD)$ be an lc pair. Find all
  Q-modifications  $\pi:S'\to S$ with reduced exceptional divisor $E'$ such that
\begin{enumerate}
\item $K_{S'}$ is $\pi$-ample, and
\item $K_{S'}+d(D'+E')\simq 0$.
\end{enumerate}

Note that the conditions
$$
K_{S'}+dD'+\Delta\simq 0\qtq{and} \supp\Delta\subset E'
\eqno{(\ref{task.task}.3)}
$$
uniquely determine $\Delta$. If  (\ref{task.task}.2) holds then
$\Delta=dE'$, and each irreducible component of $E'$ gives a linear equation for $d$. We should thus expect no solutions  if there are several irreducible components.

We start by writing down solutions of (\ref{task.task}.3), and then using this to prove that there very few other solutions.
\end{task}

\begin{say}[Discrepancy divisors] \label{P.mod.P.mod.c.lem}
  There are some easy  solutions of (\ref{task.task}.3).
\begin{enumerate}
\item
  $(S_{n,q}, D_{n,q})$.    Let $\pi:S'\to S_{n,q}$ be a $Q$-modification
  with reduced exceptional divisor $E'$. Then
  $K_{S'}+D'+E'\simq 0$.
 Indeed, this holds for the minimal resolution by (\ref{ord.of.K+C.say}.4), and continues to hold  since we blow up only nodes.
\item   $(S^d_{n,q}, D^d_{n,q})$.
  On the minimal resolution $S^{\rm m}\to S:=S^d_{n,q}$, we have
  $$
  K_{S^{\rm m}}+D^{\rm m}_{n,q}+ \tfrac12(C'_0+C''_0)+\tsum _{i=1}^s C_i\simq 0
  $$
  by (\ref{ord.of.K+C.say}.5).
  As we blow up nodes, new curves appear with  half integer coefficients that are
  $\leq 1$.
\item $(S_{n,q}, B_{n,q})$.  The formula is quite complicated already for the minimal resolution, see \cite[3.32--34]{kk-singbook}.
  We will not use these. Instead, we will choose a divisor  $\bar B_{n,q}$ so that $B_{n,q}+ \bar B_{n,q}=D_{n,q}$ and work with (\ref{P.mod.P.mod.c.lem}.1) instead.
  \end{enumerate}
\end{say}

The solution of  (\ref{task.task}) is especially simple for the pairs  $(S_{n,q}, dD_{n,q})$.

\begin{lem} \label{Dnq.triv.lem} Let $\pi:S'\to S:=S_{n,q}$ be a Q-modification
  with reduced exceptional curve $E'$ such that
  $K_{S'}$ is $\pi$-ample. Then
  \begin{enumerate}
  \item $K_{S'}+D'_{n,q}+E'\simq 0$, and
  \item if  $K_{S'}+d(D'+E')\simq 0$ for some $d<1$, then $\pi$ is an isomorphism.
    \end{enumerate}
\end{lem}

  Proof. The first claim was already noted in  (\ref{P.mod.P.mod.c.lem}.1).
  If $\pi$ has  $r\geq 1$ exceptional curves, then the  dual graph of $(S', D'+E')$ is
$$
   \bullet \ {-\boxed{\tfrac{n_0}{q_0}}-} \  e'_1 \  {-\boxed{\tfrac{n_1}{q_1}}-} \cdots  {-\boxed{\tfrac{n_{s-1}}{q_{s-1}}}-} \  e'_s \ {-\boxed{\tfrac{n_s}{q_s}}-} \ \bullet
   $$
   Since  $K_{S'}+D'+E'\simq 0$ and $K_{S'}+d(D'+E')\simq 0$, we get
   that $D'+E'\simq 0$ for  $d<1$. Then 
   $(D'+E')\cdot E'_1=0$ and  $K_{S'}\cdot E'_1=0$, a contradiction.\qed

\begin{cor} \label{34.1.no.defs.thm} A pair 
  $(S_{N,Q}, dD_{N,Q})$ is KSBA smoothable iff $S_{N,Q}$ is KSB smoothable, hence
  $N=rn^2, Q=arn-1$. If these hold then
  the universal KSBA deformation is given by
  $$
  \bigl((xy=z^{rn}+\tsum_{i=0}^{r-1}t_iz^{in}), d(z=0)\bigr)/\tfrac1{n}(1, -1, a, {\mathbf 0}).
  \hfill\qed
  $$
\end{cor}

Similarly, the dihedral pairs $(S^d_{n,q}, dD^d_{n,q})$ have no KSBA smoothings for $d<1$.  The $d=1$ case is treated in (\ref{DP.char.lem}). 

    \begin{lem} \label{Cdnq.triv.lem} Let $\pi:S'\to S:=S^d_{n,q}$ be a Q-modification such that
     $K_{S'}$ is $\pi$-ample. Assume that $K_{S'}+d(D'+E')\simq 0$ for some
     $\frac12< d<1$. Then $\pi$ is an isomorphism.
\end{lem}

     Proof.   By (\ref{P.mod.P.mod.c.lem}.2) there are half integers
     $e'_i\leq 1$ such that  $K_{S'}+D'+\sum_i e'_iE'_i\simq 0$.
     Subtracting gives that
     $$
     (1-d)D'+\tsum (e'_i-d)E'_i\simq 0.
     $$
     Since  $D'$ is $\pi$-nef, $-\tsum (e'_i-d)E'_i$ is $\pi$-nef,
     hence $e'_i-d\geq 0$ for every $i$.   Thus in fact $e'_i=1$ for every $i$.
     The rest is now the same as in (\ref{Dnq.triv.lem}), but $S^d_{n,q}$ itself has no KSB smoothings. \qed

\medskip
By contrast, the pairs  $(S_{n,q}, dB_{n,q})$ have a much more complicated behavior.
   
   \begin{lem} \label{Bnq.almosttriv.lem} Let $\pi:S'\to S:=S_{n,q}$ be a Q-modification such that
     $K_{S'}$ is $\pi$-ample. Assume that $K_{S'}+d(B'+E')\simq 0$ for some
     $ d<1$. Then $\pi$ has at most 1 exceptional curve.
\end{lem}

   Proof.  If there are at least 2 exceptional curves, then 
   the dual graph  is
$$
   \bullet \ {-\boxed{\tfrac{n_0}{q_0}}-} \  e'_1 \  {-\boxed{\tfrac{n_1}{q_1}}-} e'_2 \ - \cdots
   $$
   By (\ref{P.mod.P.mod.c.lem}.3)  $K_{S'}+B'+E'+\bar B'\simq 0$.
   Subtracting, we get that
   $(1-d)(B'+E')+\bar B'\simq 0$.
   Here $\bar B'$ is disjoint from $E'_1$, so
   $(B'+E')\cdot E'_1=0$.  This in turn gives that 
$K_{S'}\cdot E'_1=0$, a contradiction.\qed

\begin{exmps} \label{4.0.smooth.exmp}
  Note that $S'$ with dual graph 
$$
   \bullet \ {-\boxed{\tfrac{n_0}{q_0}}-} \  e'_1 \  {-\boxed{\tfrac{n_1}{q_1}}}
   $$
  has a KSB deformation with Du~Val general fiber iff $S_{n_0, q_0}$ is a T-singularity and
   $S_{n_1, q_1}$ is a  Du~Val singularity. 
The simplest case  is the dual graph
     $$
   \bullet\ {-\boxed{\tfrac{4}{1}}-} \  c. 
    $$
   This gives that
   $
   \bigl(\a^2_{xy},\tfrac{2c-3}{2c-1}(y=0)\bigr) /\tfrac{1}{4c-1}(1, c)
     $
   is KSBA smoothable.

   I found only 2 series of pairs whith KSBA deformations for different values of $d$. The first is
   $ \bullet\ - \ 4\ - \  3  \ {-\boxed{A_{n+2}}}$
with P-modifications 
   $$
   \bullet\ {-\boxed{\tfrac{4}{1}}-} \  3  \ {-\boxed{A_{n+2}}}
   \qtq{and} 
   \bullet\ {-\boxed{\tfrac{18}{5}}-} \  2  \ {-\boxed{A_{n}}}.
   $$
  The second is  $ \bullet\ - \ 2\ - \ 5\ - \  3  \ {-\boxed{A_{n+1}}}$
with P-modifications 
   $$
   \bullet\ {-\boxed{\tfrac{9}{5}}-} \  3  \ {-\boxed{A_{n+1}}}
   \qtq{and} 
   \bullet\ {-\boxed{\tfrac{25}{14}}-} \  2  \ {-\boxed{A_{n}}}.
   $$
   \end{exmps}

\begin{exmp} \label{D.term.QC.exmp.1} Doubly KSB deformations of $D_n$ singularities are described as follows.
  For $r\geq 2$ the $D_r$-type
  pair
  $\bigl((z^2=x(y^2-x^{r-2})), (x=z=0)\bigr)$ has an 
$(r{-}2)$-parameter deformation 
  $$
  z^2=x\bigl(y^2-x^{r-2}-\tsum_{i=0}^{r-3} t_i x^i\bigr).
  $$
  (For $r=2$ we mean $(z^2=x(y^2-1))$, which has 2 singular points of type
  $\bigl((z^2=xy), (x=z=0)\bigr)$ at $(0, \pm 1, 0)$.)
  
  The divisor  $D:=(x=z=0)$ is $\q$-Cartier, since $2D=(x=0)$.
  So these are doubly KSB deformations.
The generic deformation  has 2 singular points of type
  $\bigl((z^2=xy), (x=z=0)\bigr)$ at   $(0,\pm \sqrt{-t_0}, 0)$.
  By (\ref{ord.of.K+C.say}.1) the latter has no KSB smoothings.
  \end{exmp}

\begin{exmp} \label{1/2.2.exmp} There are a several   KSBA smoothable
  examples with $d=\frac12$. These are given by the dual graph, where $\ast$ can be either $\bullet$ or $2$
    $$
  \begin{array}{ccccc}
    && \ast &&\\
&& | &&\\
    \ast & - & 4  & - &   \ast\\
    && | &&\\
     && \ast &&
  \end{array}
  \qtq{where $D^{\rm m}+\Delta^{\rm m}$ is}
  \begin{array}{ccccc}
    && \tfrac12 &&\\
&& | &&\\
     \tfrac12 & - & 1  & - &    \tfrac12\\
    && | &&\\
     &&  \tfrac12 &&
  \end{array}
  $$
  Contracting the curve $(4)$ we get 
   doubly KSBA smoothable examples.
\end{exmp}

\section{First order deformations}\label{first.ord.sec}

Various aspects of the first order deformation theory of cyclic quotients are worked out in \cite{k-alt}; the references below are to the version in
\cite[Sec.6.6]{k-modbook}.

Fix a pair  $(S, D):=\bigl(\a^2_{xy}, (xy=0)\bigr)/\tfrac1{n}(1,q)$.
Using  $\dx:=x(\partial/\partial x)$ and  $\dy:=y(\partial/\partial y)$,
the first order deformation space  $T^1(S)$ of $S$ has a basis consisting of certain operators
$$
{\mathcal D}=\tfrac{\alpha\dx-\beta\dy}{M},\qtq{where} \alpha,\beta\in \c\qtq{and}
M=x^ay^b.
$$
Since  $\omega_S(D)\cong \o_S$, we can think of the curve $D$ as determined by
$$
\o_S\to \omega_S^{-1}\qtq{given by}  1\mapsto \tfrac{xy}{dx\wedge dy}.
$$
The lifting criterion \cite[6.58]{k-modbook} says that
$\frac{xy}{dx\wedge dy} $ lifts to a first order deformation determined by
 ${\mathcal D}=({\alpha\dx-\beta\dy})/(x^ay^b)$ iff
$$
{\mathcal D}(xy)-xy\nabla {\mathcal D}= (b\beta-a\alpha) \tfrac{xy}{x^ay^b}
$$
equals $0$, as an element in $H^1(\a^2\setminus\{(0,0)\}, \o)$.
The latter holds if either $b\beta-a\alpha=0$ or $\min\{a,b\}\leq 1$.

Once $\frac{xy}{dx\wedge dy} $ lifts to a first order deformation
$S[\epsilon]$, the lift $D[\epsilon]$ of $D$  is given by the induced
$\omega_{S[\epsilon]/k[\epsilon]}\to \o_{S[\epsilon]}$, thus
$D[\epsilon]$ is also flat over $k[\epsilon]$. Hence
$\bigl(S[\epsilon], D[\epsilon]\bigr)$ is a KSB deformation of $(S, D)$ as in  \cite[6.21--22]{k-modbook}.

We can now look at the basis of  $T^1(S)$ given in \cite[6.73]{k-modbook}.
The above condition is satisfied by the basis elements $\dx/M_1, \dy/M_{r-1}$,
by a 1-dimensional subspace  of each  $\bigl\langle \dx/M_i, \dy/M_i\bigr\rangle$ for $2\leq i\leq r{-}2$, and by all the basis elements in  the list  \cite[6.73.2]{k-modbook}. Since the $r$ in \cite[6.73]{k-modbook} is $m(S)-1$ (where $m(S)$ is the multiplicity), we have proved the following.

\begin{thm}\label{KSB.codin.infini.thm}
  Set $(S, D):=\bigl(\a^2_{xy}, (xy=0)\bigr)/\tfrac1{n}(1,q)$. Then $T^1_{\rm KSB}(S, D)$ has codimension $m(S)-3$ in
  $T^1(S)$. \qed
  \end{thm}

\begin{exmp} We compute in more detail what happens for the singularity
  $S_{n,1}:=\a^2/\frac1{n}(1,1)$ for $n\geq 3$.

  The end result is that $\dim T^1(S_{n,1})=2n-4$, the Artin component is smooth
  of dimension $n-1$, and $T^1_{\rm KSB}(S_{n,1}, D_{n,1})$  also
has dimension $n-1$. Together, $T^1_{\rm Art}(S_{n,1})$ and $T^1_{\rm KSB}(S_{n,1}, D_{n,1})$
span $T^1(S_{n,1})$, so their intersection has dimension 2.

By Theorem~\ref{main.cyclic.thm}, $\defsp_{\rm KSB}(S_{n,1}, D_{n,1})$ has dimension 2; it is explicitly worked out  in   \cite[2.32]{k-modbook}.
In particular, $\defsp_{\rm KSB}(S_{n,1}, D_{n,1})$ is nonreduced for $n\geq 4$.

To see these, note first  that by
\cite{pinkham}, the Artin component can be given by the  equations 
$$
\rank
\left(
\begin{array}{ccccc}
  z_0 & z_1 & \cdots &  z_{n-2} & z_{n-1}\\ 
 z_1-t_1 & z_2-t_2 & \cdots &  z_{n-1}-t_{n-1} & z_{n}
\end{array}
\right) \leq 1.
$$
The coordinate ring of $S_{n,1}$ is generated by the monomials
$M_i:=x^{n-i}y^i$. We claim that $T^1_{\rm Art}(S_n)$ is spanned
 by the vector fields
$$
{\mathcal D}_i:=\tfrac{(n-i)\dx+i\dy}{nM_i}.
$$
We compute that 
${\mathcal D}_i(M_j)=M_j/M_i$,
which is invertible on the $x\neq 0$ chart for $j\geq i$. Thus, in this chart, we can choose 
$$
M_0, M_1, \dots, M_{i-1}, M_i+\epsilon Q_i,
M_{i+1}+\epsilon Q_{i+1},\dots,
M_n+\epsilon Q_n
$$
as liftings of the $M_i$, where we set  $Q_j:={M_j}/{M_i}$.
Note that  the matrix
$$
\left(
\begin{array}{ccccccccc}
  M_0 & M_1 & \cdots & M_{i-1} & M_i+\epsilon Q_i & \cdots & M_{n-2}+\epsilon Q_{n-2} & M_{n-1}+\epsilon Q_{n-1}\\
  M_1 & M_2 & \cdots & M_i & M_{i+1}+\epsilon Q_{i+1}&\cdots &  M_{n-1}+\epsilon Q_{n-1} & M_n+\epsilon Q_n
\end{array}
\right)
$$
has rank $\leq 1$. Relabeling the entries of the top row as
$z_0,\dots, z_{n-1}$, these become  the equations
$$
\rank
\left(
\begin{array}{ccccccccc}
  z_0 & z_1 & \cdots & z_{i-1} & z_i & \cdots &  z_{n-2} & z_{n-1}\\ 
 z_1 & z_2 & \cdots & z_{i}-\epsilon & z_{i+1} & \cdots &  z_{n-1} & z_{n}
\end{array}
\right) \leq 1.
$$
Thus these deformations span the Artin component.
\end{exmp}

The philosophy of \cite{k-modbook} says that in order to get the optimal moduli theory, we should  establish properties of 1-parameter deformations, and then formulate versions of these   properties  over arbitrary bases.
The second step is not necessarily unique. 
For deformations of pairs, this poses the following.

\begin{ques} Is there a good way to impose additional requirements for KSB deformations of a pair to make  $\defsp_{\rm KSB}(S_{n,q}, D_{n,q})$ smooth?
    \end{ques}

\section{Deformations of modifications}\label{def.mod.sect}

We prove that Q-modifications $X\to S$ have no (small or large)  KSB deformations over reduced bases that keep $S$ fixed. (See (\ref{An.mines.exmp}) for other deformations and non-reduced bases.)
More generally, we have the following.

\begin{thm} \label{rigid.quot.mod.thm}
  Let $S$ be a normal surface, $B$ a smooth, irreducible curve, and $g:Y\to S\times B$ a proper birational morphism, $Y$ normal. Assume that
  \begin{enumerate}
  \item $g_b:Y_b\to S$ is birational for every $b\in B$, and
    \item there is a dense,   open subset  $B^\circ\subset B$ such that, for every $b\in B^\circ$,
  $Y_b\to S$  is a Q-modification.
      \end{enumerate}
  Then, $Y\to B$ is trivial. That is, for any $b\in B$, 
    $$
  \bigl(Y\to S\times B\to B\bigr)\cong
  \bigl(Y_b\times B\to S\times B\to B\bigr).
  \eqno{(\ref{rigid.quot.mod.thm}.3)}
  $$
\end{thm}

Proof.   Over a dense,  open subset $B^*\subset B^\circ$,
the minimal resolutions   $\{Y_b^{\rm m}\to Y_b:b\in  B^*\}$
form a flat family. By  (\ref{rigid.quot.mod.lem.1}),
we have only finitely many choices for the centers of the blow-ups giving
$Y_b^{\rm m}\to S^{\rm m}$, so the $\{Y_b^{\rm m}:b\in  B^*\}$,
and hence also $\{ Y_b:b\in  B^*\}$, form localy trivial families.

Thus the family is also trivial over $B$
by (\ref{rigid.quot.mod.lem.2}). \qed

\begin{lem} \label{rigid.quot.mod.lem.2}
  Let $g:Y\to X\times B$ be a projective morphism such that
  $Y$ is normal and 
  $g_b:\red (Y_b)\to X_b$ is birational for every $b$.
  Assume that there is a dense, open $B^\circ\subset B$ such that
  $(Y^\circ\to B^\circ)\cong ((Y'\times B^\circ)\to B^\circ)$ for some
  normal $Y'$. Then $(Y\to B)\cong ((Y'\times B)\to B)$.
\end{lem}

Proof. Let $H$ be relatively ample on $Y\to X\times B$ and let $H'$ be its
birational  transform on $(Y'\times B)\to B$.
$H'$ is Cartier over $B^\circ$, so Cartier by \cite{ram-sam, sam-ram}; see also
\cite[4.21]{k-modbook}.
Thus $H'$ is relatively ample on $(Y'\times B)\to B$.
Thus $(Y\to B)\cong ((Y'\times B)\to B)$
 by   \cite{mats-mumf}; see \cite[11.39]{k-modbook} for the form that we use.  \qed

 \begin{exmp}\label{An.mines.exmp} M-modifications have nontrivial deformations, as
   the next examples show.

   (\ref{An.mines.exmp}.1)  Let $S_P\to S$ be any P-modification whose construction involves blowing up a point  $p\in S^{\rm m}$. We can move the point
   $p$ along the exceptional curve and contract in the family to get a
   nontrivial flat deformation of $S_P\to S$ (keeping $S$ fixed).
   Here the canonical class of the generic fiber is not relatively nef.
   The simplest example is  the singularity $3\ -\ 3$. After one blow-up, the central and generic fibers are
   $$
   4\ -\ 1\ -\ 4 \qtq{resp.}   3\ -\ 4\ -\ 1.
   $$
   All curves marked $4$ are then contracted. The central fiber is an 
   M-modification. In the generic fiber, the canonical class has degree $-\frac12$ on the image of the
   curve marked $1$. 
   
    (\ref{An.mines.exmp}.2) For $A_n:=(xy=z^{n+1})$, the miniversal deformation of $A_n$ is
   $$
   \bigl(xy=z^{n+1}+\tsum_{i=0}^{n-1}t_iz^i\bigr),
   $$
   while the  miniversal deformation  of its minimal  resolution  is  obtained from 
   $$
   \bigl(xy=\tprod_{i=0}^{n}(z-s_i)\bigr), \qtq{(subject to $\tsum s_i=0$)}
   $$
   by  repeatedly blowing up the divisors $(x=z-s_i=0)$ as in (\ref{ksb.defs.say}.3).
   
   Let $\sigma_j$ be the $j$th elementary symmetric polynomial in the $s_i$.
   Then, over the Artin ring  $k[s_0,\dots, s_n]/(\sigma_0, \dots, \sigma_n)$, we get
   a   deformation  of the minimal  resolution which contracts to the
   trivial deformation of  $A_n$.

   A similar result holds for all  Du~Val singularities by \cite{art-bri},
   and for M-resolutions by \cite{be-ch}.
   \end{exmp}


 \def\cprime{$'$} \def\cprime{$'$} \def\cprime{$'$} \def\cprime{$'$}
  \def\cprime{$'$} \def\dbar{\leavevmode\hbox to 0pt{\hskip.2ex
  \accent"16\hss}d} \def\cprime{$'$} \def\cprime{$'$}
  \def\polhk#1{\setbox0=\hbox{#1}{\ooalign{\hidewidth
  \lower1.5ex\hbox{`}\hidewidth\crcr\unhbox0}}} \def\cprime{$'$}
  \def\cprime{$'$} \def\cprime{$'$} \def\cprime{$'$}
  \def\polhk#1{\setbox0=\hbox{#1}{\ooalign{\hidewidth
  \lower1.5ex\hbox{`}\hidewidth\crcr\unhbox0}}} \def\cdprime{$''$}
  \def\cprime{$'$} \def\cprime{$'$} \def\cprime{$'$} \def\cprime{$'$}
\providecommand{\bysame}{\leavevmode\hbox to3em{\hrulefill}\thinspace}
\providecommand{\MR}{\relax\ifhmode\unskip\space\fi MR }
\providecommand{\MRhref}[2]{%
  \href{http://www.ams.org/mathscinet-getitem?mr=#1}{#2}
}
\providecommand{\href}[2]{#2}

 \bigskip

  Princeton University, Princeton NJ 08544-1000, \

  \email{kollar@math.princeton.edu}

 \end{document}